\documentclass[12pt]{amsart}

\usepackage{amscd,latexsym,amsthm,amsfonts,amssymb,amsmath,amsxtra}
\usepackage[all]{xy}

\usepackage[mathscr]{eucal}
\usepackage[pdftex,bookmarks,colorlinks,pdfmenubar]{hyperref}

\renewcommand\theequation{\thesection.\arabic{equation}}

\newcommand{\BA}{{\mathbb {A}}}

\newcommand{\BC}{{\mathbb {C}}}

\newcommand{\CA}{{\mathcal {A}}}
\newcommand{\CB}{{\mathcal {B}}}
\newcommand{\CC}{{\mathcal {C}}}

\newcommand{\CE}{{\mathcal {E}}}
\newcommand{\CF}{{\mathcal {F}}}
\newcommand{\CG}{{\mathcal {G}}}
\newcommand{\CH}{{\mathcal {H}}}

\newcommand{\CO}{{\mathcal {O}}}

\newcommand{\CZ}{{\mathcal {Z}}}

\newcommand{\Fd}{{\mathfrak {d}}}

\newcommand{\Fg}{{\mathfrak {g}}}

\newcommand{\Fk}{{\mathfrak {k}}}
\newcommand{\Fl}{{\mathfrak {l}}}
\newcommand{\Fm}{{\mathfrak {m}}}
\newcommand{\Fn}{{\mathfrak {n}}}

\newcommand{\Fq}{{\mathfrak {q}}}
\newcommand{\Fr}{{\mathfrak {r}}}

\newcommand{\RD}{{\mathrm {D}}}

\newcommand{\RO}{{\mathrm {O}}}

\newcommand{\Ad}{{\mathrm{Ad}}}
\newcommand{\Aut}{{\mathrm{Aut}}}

\newcommand{\cusp}{{\mathrm{cusp}}}

\newcommand{\disc}{{\mathrm{disc}}}

\newcommand{\GL}{{\mathrm{GL}}}

\newcommand{\Hom}{{\mathrm{Hom}}}

\newcommand{\Int}{{\mathrm{Int}}}

\newcommand{\out}{{\mathrm{Out}}}

\newcommand{\PGL}{{\mathrm{PGL}}}
\newcommand{\PD}{{\mathrm{PD}}}

\renewcommand{\Re}{{\mathrm{Re}}}

\newcommand{\simp}{{\mathrm{sim}}}

\newcommand{\SO}{{\mathrm{SO}}}

\newcommand{\Sp}{{\mathrm{Sp}}}

\newcommand{\st}{{\mathrm{st}}}
\newcommand{\Span}{{\mathrm{Span}}}

\newcommand{\ud}{\,\mathrm{d}}

\newcommand{\udl}{\underline}
\newcommand{\wt}{\widetilde}

\newcommand{\cpair}[1]{\left\{{#1}\right\}}

\newcommand{\ol}{\overline}

\newcommand{\bs}{\backslash}

\def\bks{{\backslash}}

\def\diag{{\rm diag}}

\def\sym{{\rm sym}}
\def\sig{{\sigma}}

\newtheorem{thm}{Theorem}[section]

\newtheorem{prop}[thm]{Proposition}
\newtheorem {conj}[thm]{Conjecture}
\newtheorem {ques/conj}[thm]{Question/Conjecture}
\newtheorem {ques}[thm]{Question}

\newtheorem{prob}[thm]{Problem}

\usepackage[numbers]{natbib}

\makeatletter

\newcommand{\Rmnum}[1]{\expandafter\@slowromancap\romannumeral #1@}
\makeatother

\begin{document}
\renewcommand{\theequation}{\arabic{equation}}
\numberwithin{equation}{section}

\title[Bessel Descents and Branching Problems]{Bessel Descents and Branching Problems}

\author{Dihua Jiang}
\address{School of Mathematics, University of Minnesota, Minneapolis, MN 55455, USA}
\email{dhjiang@math.umn.edu}

\author{Lei Zhang}
\address{Department of Mathematics,
National University of Singapore,
Singapore 119076}
\email{matzhlei@nus.edu.sg}

\subjclass[2000]{Primary 11F67, 11F70, 22E55; Secondary 11F30, 11F66}

\keywords{Endoscopic Classification, Discrete Spectrum, Bessel Periods, Twisted Automorphic Descents, Branching Problem}

\date{\today}
\thanks{The research of the first named author is supported in part by the NSF Grant DMS--1600685 and DMS-1901802, and that 
of the second named author is supported in part by AcRF Tier 1 grant R-146-000-237-114 and R-146-000-277-114 of National University of Singapore. }

\begin{abstract}
We discuss the theory of automorphic descents of Bessel type and its relation to automorphic version of branching problem and its relevant reciprocal branching problem. 
\end{abstract}

\maketitle


\section{Introduction}


Automorphic descent is a method to construct automorphic forms on a smaller group from these on a larger group. The automorphic descent of Bessel type is a construction of such descents using the Bessel-Fourier 
coefficients of automorphic forms on a larger group, which is simply called the {\sl Bessel descent} as in the title of this paper. There is another type of the automorphic descent, which may be called 
{\sl Fourier-Jacobi descent}, using Fourier-Jacobi coefficients of automorphic forms. 
The idea goes back to the construction of the Saito-Kurokawa cusp forms of 
$\Sp_4$, which are non-tempered at almost all local places and hence are counter-examples to the generalized Ramanujan conjecture (\cite{An79}, \cite{PS83}, and \cite{EZ85}). 

The more recent development of 
the automorphic descent was re-initiated by the work of D. Ginzburg, S. Rallis and D. Soudry in the 1990's in order to construct the backwards global Langlands functorial transfers 
from the general linear groups to quasi-split classical groups. The idea is to construct explicit cuspidal automorphic modules on classical groups by means of their global Arthur parameters. 
Their work has been well presented in their book (\cite{GRS11}). Note that the automorphic descents in \cite{GRS11} requires that the classical groups are quasi-split and 
the cuspidal automorphic representations on quasi-split classical groups must be generic in the sense that they admit a non-zero Whittaker-Fourier coefficients. 

The notion of twisted automorphic descents was introduced by the authors in their previous work (\cite{JZ-BF}, and also \cite{JZ-Howe}), which can be viewed as a natural extension of the automorphic descents of Ginzburg-Rallis-Soudry. The twisted automorphic descent is to take, as source representations, square-integrable automorphic representations $\pi$ of classical group $G_n$ with global Arthur parameters of type 
$$
\psi^{(2,1)}=\phi_\tau^{(2)}\boxplus\phi_\delta.
$$
Here $\phi_\tau$ and $\phi_\delta$ are generic global Arthur parameters of the form in \eqref{aps} with all $b_j=1$, and $\phi_\tau^{(2)}$ is defined to be the following non-generic global Arthur parameter:
\begin{equation}\label{aps2}
\phi_\tau^{(2)}=(\tau_1,2)\boxplus\cdots\boxplus(\tau_r,2),
\end{equation}
if $\phi_\tau=(\tau_1,1)\boxplus\cdots\boxplus(\tau_r,1)$ with $\tau=\tau_1\boxplus\cdots\boxplus\tau_r$. We refer to Section \ref{ssec-ds} for details on discrete spectrum of classical groups and their 
Arthur parameters. 

The square-integrable residual representations of quasi-split classical group $G_n^*$ with global Arthur parameter 
$\psi^{(2,1)}$ are explicitly determined in \cite{JLZ13}. Note that the automorphic descent of Ginzburg-Rallis-Soudry (\cite{GRS11}) only considers the residual representations with global Arthur parameters of the form 
$\phi_\tau^{(2)}$ and only for quasi-split classical groups. The key point in the theory of twisted automorphic descents is to introduce a generic global Arthur parameter $\phi_\delta$, with the property that 
the central $L$-value $L(\frac{1}{2},\tau\times\delta)$ of the twisted $L$-functions is non-zero. 
This property plays an indispensable role in the construction of the twisted automorphic descents for general classical groups. 

The twisted automorphic descents have been proved to be fruitful and productive in the endoscopy theory of automorphic forms and its applications. 
We list the following important development in the theory of twisted automorphic descents:
\begin{enumerate}
\item In \cite{JLXZ}, a special case of the twisted automorphic descents was explicitly studied, which leads to a new proof of the classical Jacquet-Langlands correspondence between $\PGL_2$ and $\PD^\times$, 
where $\RD$ is a quaternion algebra defined over a number field $k$. 
\item In \cite{JZ-BF}, the twisted automorphic descents were systematically developed in order to construct explicit cuspidal automorphic modules for general classical groups by means of their global Arthur 
parameters, with assumption that the global Arthur parameters are generic, and with the assumption of the {\sl Generic Summand Conjecture} (Conjecture \ref{conj-gs}). 
\item In \cite{JZ-BF}, the twisted automorphic descents (ideas and methods) are used to proved  one direction of the global Gan-Gross-Prasad conjecture unitary groups and special orthogonal groups in a uniform way, and 
prove the other direction of the global Gan-Gross-Prasad conjecture with a global assumption that can only be verified in some special cases, currently. 
\item In \cite{JZ18}, the idea and method of the twisted automorphic descents was used to prove the non-vanishing of certain tensor product $L$-functions for classical groups. 
\item In \cite{JLX}, the twisted automorphic descents were used to test the {\sl Branching Problem} and the {\sl Reciprocal Branching Problem} for cuspidal automorphic representations of classical groups. 
\end{enumerate}
We will discuss in this paper the theory of twisted automorphic descents and its relation to the automorphic version of the branching problem and its relevant reciprocal branching problem. 
Refer to \cite{J14}, \cite{J-Shahidi}, and \cite{JZ-Howe} for further discussions on other relevant topics. 
Instead of treating the theory for general classical groups, we mainly discuss the theory of Bessel descents 
for special orthogonal groups and leaves more detailed treatment of other classical groups in our future work. 

The branching decomposition problem roots in the classical theory of representations and invariants. 
Let $G$ be a group and $H$ be a subgroup of $G$. The classical {\sl Branching Decomposition Problem} is to ask: For a given representation $\pi$ of $G$, what is the spectrum of the restriction of $\pi$ to $H$ 
as a representation of $H$? The {\sl Reciprocal Branching Problem} is then to ask: For a given representation $\sigma$ of $H$, how to find a representation $\pi$ of $G$, such that the restriction of $\pi$ to $H$ 
contains $\sigma$? It is clear that without any extra conditions or requirements on $\pi$ and/or $\sigma$, the two problems are related by the well known {\sl Frobenius Reciprocity}. 
However, with some extra requirements on $\pi$ and $\sigma$, it could become an interesting and possibly difficult problem. 

Let $G_n$ be a special orthogonal group over a number field $k$ as defined in Section \ref{sec-BMBP}, although the discussion here makes sense for general classical groups, or even for general reductive algebraic groups. 
Let $H_m$ be a special orthogonal group over a number field $k$, which has a natural embedding into $G_n$ as a subgroup, so that the pair $(G_n,H_m)$ is relevant in the sense of the global Gan-Gross-Prasad conjecture 
(\cite{GGP12}).
For any irreducible unitary representation $\pi$ of $G_n(\BA)$ (where $\BA$ is the ring of adeles of $k$), which has a realization $\CC_\pi$ in the space of the automorphic discrete spectrum of $G_n$ 
(see Section \ref{ssec-ds} for definition), 
the automorphic version of {\sl branching problem} may ask: what is the spectral decomposition of the space $\CC_\pi$ restricted to $H_m(\BA)$ as an automorphic 
$H_m(\BA)$-module? It is in general a complicated problem if $H_m$ is small comparing to $G_n$, since the decomposition may not be discrete, and the multiplicity may be infinite, for instance. 
From representation-theoretic point of view, it will be reasonable to consider first the case of multiplicity one situation. To this end, one may enlarge the subgroup $H_m$ in a reasonable way. 

One of the natural ways is to consider the centralizer $C_{G_n}(H_m)$ of the subgroup $H_m$ in $G_n$, which is a reductive subgroup of $G_n$. It makes sense to consider the natural embedding:
$$
C_{G_n}(H_m)\times H_m\rightarrow G_n.
$$
It is usually expected that one may consider the branching decomposition of the space $\CC_\pi$ restricted to $C_{G_n}(H_m)(\BA)\times H_m(\BA)$ as an intermediate step towards the eventual understanding of the branching 
decomposition problem from $G_n$ to the subgroup $H_m$. This leads to an interesting topic, but beyond the reach of this paper. 

In this paper, we consider another way to enlarge the subgroup $H_m$. Namely, we consider,  
when $m=\ell^-$, the Bessel subgroup of $G_n$
\begin{equation}\label{bsg}
R_\ell:=H_{\ell^-}^{\CO_\ell}\rtimes V_\ell,
\end{equation}
where $H_m:=H_{\ell^-}^{\CO_\ell}$ and a unipotent subgroup $V_\ell$ of $G_n$ are defined in Section \ref{ssec-bmbp}. The Bessel-Fourier coefficient of $\CC_\pi$ attached to the Bessel subgroup produces 
an automorphic Bessel module $\CF^{\CO_\ell}(\CC_\pi)$ of $H_m(\BA)$ as defined in \eqref{bm}. The question is: how to understand the spectral structure of the automorphic Bessel module $\CF^{\CO_\ell}(\CC_\pi)$ 
as representation of $H_m(\BA)$? In particular, we are interested in the following branching problem. 
\begin{prob}[$L^2$-Automorphic Branching Problem]\label{q-L2bp}
Let $\pi$, $\CC_\pi$, $G_n$, $H_m$, and $\CF^{\CO_\ell}(\CC_\pi)$ be as given above. 
Find irreducible unitary representations $\sigma$ of $H_m(\BA)$, which has a realization $\CC_\sigma$ in the space of discrete spectrum of $H_m(\BA)$, so that the Bessel period 
$\CB^{\CO_{\ell}}(\varphi_\pi,\varphi_\sigma)$, as defined in \eqref{bp-L2}, is non-zero for some $\varphi_\pi\in\CC_\pi$ and $\varphi_\sigma\in\CC_\sigma$.
\end{prob}

Note that the Bessel period $\CB^{\CO_{\ell}}(\varphi_\pi,\varphi_\sigma)$, if it is non-zero, defines a non-zero linear functional in 
$$
\Hom_{H_{\ell^-}^{\CO_{\ell}}(\BA)}(\CF^{\CO_{\ell}}(\CC_\pi)\otimes\sigma^\vee,1).
$$
From the uniqueness of local Bessel models (\cite{AGRS}, \cite{SZ},\cite{GGP12}, and \cite{JSZ}), for any $\sigma\in\CA_\cusp(H_{\ell^-}^{\CO_{\ell}})$, the $\Hom$-space 
is at most one-dimensional. Hence if the Bessel period $\CB^{\CO_{\ell}}(\varphi_\pi,\varphi_\sigma)$ is non-zero, $\sigma$ occurs in the spectral decomposition of the automorphic Bessel module 
$\CF^{\CO_{\ell}}(\CC_\pi)$ with multiplicity one. 

The {\sl Reciprocal Branching Problem} can be formulated as follows. 
\begin{prob}[$L^2$-Automorphic Reciprocal Branching Problem]\label{q-L2rbp}
Let $H_m$ be a special orthogonal group defined over a number field $k$. Given an irreducible unitary representation $\sigma$ of $H_m(\BA)$, 
which has a realization $\CC_\sigma$ in the space of discrete spectrum of $H_m(\BA)$, find an special orthogonal group $G_n$ containing $H_m$ as a subgroup such that $(G_n,H_m)$ is a relevant pair in the sense of 
the global Gan-Gross-Prasad conjecture, and find an irreducible unitary representation $\pi$ of $G_n(\BA)$, which has a realization $\CC_\pi$ in the space of discrete spectrum of $G_m(\BA)$, such that 
the Bessel period $\CB^{\CO_{\ell}}(\varphi_\pi,\varphi_\sigma)$, as defined in \eqref{bp-L2}, is non-zero for some $\varphi_\pi\in\CC_\pi$ and $\varphi_\sigma\in\CC_\sigma$.
\end{prob}

The objective of this paper is to understand these two problems by means of the theory of twisted automorphic descents, as developed in \cite{JZ-BF}. We expect the discussion in this paper will be applicable to other 
classical groups, according to our work \cite{J14}, \cite{JZ-BF}, \cite{J-Shahidi}, and \cite{JZ-Howe}. 

The paper is organized as follows. We recall in Section \ref{sec-BMBP}the basic statement of endoscopy classification of the discrete spectrum for special orthogonal groups (\cite{A13}), the basic definition of 
Bessel-Fourier coefficients and Bessel periods of automorphic forms on special orthogonal groups (\cite{JZ-BF}), and a family of global zeta integrals that built from Bessel periods and the basic identity of 
such global zeta integral with tensor product $L$-functions (Theorem \ref{thm-bi}). By applying Theorem \ref{thm-bi} in various situations, we are able to control the basic structure of global Arthur parameters 
of $\pi$ and $\sigma$ when they have a non-zero Bessel periods in Section \ref{sec-BPGAP}. Due to technical issues as explained in Section \ref{sec-BPGAP}, we are able to prove Theorems \ref{thm-psim}, \ref{thm-psin} 
when one of $\pi$ and $\sigma$ has a generic global Arthur parameter. Based on these results, we go further to prove Theorem \ref{thm-bpclv} that the Bessel period of $\pi$ and $\sigma$ is non-zero implies 
the central value ($s=\frac{1}{2}$) of a quotient of relevant $L$-functions is non-zero. According to the notes of Gan-Gross-Prasad (\cite{GGP18}), this can be viewed as one direction of the global 
Gan-Gross-Prasad conjecture for non-generic global Arthur parameters, although the precise formulation of such conjecture for general situation is still not known as far as the authurs know. 
As explained in \cite{JZ-BF}, the theory of twisted automorphic descents is naturally connected to the global Gan-Gross-Prasad conjecture for generic global Arthur parameters. Hence as expected, 
the further refinement of the theory of twisted automorphic descents will be a powerful method to treat the global Gan-Gross-Prasad conjecture in general. After recalling the theory of twisted automorphic descents 
in Section \ref{sec-TTAD}, we go on in Section \ref{sec-BPC} the automorphic version of the {\sl Branching Problem} for cuspidal automorphic representations with generic global Arthur parameters. This can be viewed as 
reformulation of some main results in \cite{JZ-BF} on construction of cuspidal automorphic modules for special orthogonal groups via the twisted automorphic descents. In the last section (Section \ref{sec-RBP}), 
we consider the automorphic version of the {\sl Reciprocal Branching Problem} for cuspidal automorphic representations with a given generic global Arthur parameters (Theorems \ref{thm-crb} and \ref{thm-crb-foc}). 
At the end of Section \ref{sec-RBP}, the connection of the results here with the global Gan-Gross-Prasad conjecture in general is briefly discussed. 

The research work of this paper started during the visit in January, 2018, of the first named authur to National University of Singapore. He would like to thank Math. Department for providing productive research 
environment and hospitality. 
The first named author would like to thank the organizers: W. Mueller, S. W. Shin and N. Templier for invitation to the wonderful Simons Symposium 2018 on Relative Trace Formulas and thank the Simons Foundation for hospitality.


\section{Bessel Modules and Bessel Periods}\label{sec-BMBP}


Let $k$ be a number field and $\BA$ the ring of adeles of $k$.
As in \cite{JZ-BF} and \cite[Chapter 9]{A13}, we may use $G_n^*=\SO(V^*,q^*_{V^*})$ for a $k$-quasisplit special orthogonal group that is defined by a non-degenerate, 
$\Fn$-dimensional quadratic space $(V^*,q^*)$ over $k$ with $n=[\frac{\Fn}{2}]$ and use $G_n=\SO(V,q)$ to denote a pure inner $k$-form of $G_n^*$. This means that both quadratic spaces $(V^*,q^*)$ and $(V,q)$ have the same dimension and the same discriminant, as discussed in \cite{GGP12} and \cite{JZ-BF}, for instance.

Let $(V_0,q)$ be the $k$-anisotropic kernel of $(V,q)$ with dimension $\Fd_0=\Fn-2\Fr$, where the $k$-rank $\Fr=\Fr_\Fn=\Fr(G_n)$ of $G_n$ is
equal to the Witt index of $(V,q)$.
Let $V^{+}$  be a maximal totally isotropic subspace of $(V,q)$, with a basis $\{e_{1},\dots,e_{\Fr}\}$. 
Then $(V,q)$ has the following polar decomposition
$$
V=V^{+}\oplus V_{0}\oplus V^{-},
$$
where $V_{0}=(V^{+}\oplus V^{-})^{\perp}$. Let $\{e_{-1},\dots,e_{-\Fr}\}$ be a basis of $V^{-}$ that is dual to the basis $\{e_{1},\dots,e_{\Fr}\}$, i.e. 
$
q(e_{i},e_{-j})=\delta_{i,j}
$
for all $1\leq i,j\leq\Fr$. We may choose an orthogonal basis $\{e'_{1},\dots,e'_{\Fd_0}\}$ of $V_{0}$ with the property that
$
q(e'_{i},e'_{i})=d_{i},
$
where $d_{i}$ is nonzero for all $1\leq i\leq \Fd_0$. Hence we obtain the following basis for $(V,q)$:
\begin{equation}\label{bs}
e_{1},\dots,e_{\Fr},e'_{1},\dots,e'_{\Fd_0},e_{-\Fr},\dots,e_{-1}.
\end{equation}
We may fix a full isotropic flag in $(V,q)$:
$$
\Span\{e_{1}\}\subset\Span\{e_{1},e_{2}\}\subset
\cdots\subset
\Span\{e_{1},\dots,e_{\Fr}\},
$$
which defines a minimal parabolic $k$-subgroup $P_0$.
Moreover, $P_0$ contains a maximal $k$-split torus $S$, consisting of elements
$$
\diag\{t_{1},\dots,t_{\Fr},1,\dots,1,{t}^{-1}_{\Fr},\dots,t^{-1}_{1}\},
$$
with $t_i\in k^\times$ for $i=1,2,\cdots,\Fr$. Then the centralizer $Z(S)$ in $G_n$ is $S\times G_{d_0}$, the Levi subgroup
of $P_0$, where $G_{d_0}$ is defined by the $k$-anisotropic $(V_0,q_0)$, which has dimension $\Fd_0$, and $d_0=[\frac{\Fd_0}{2}]$. 
Then $P_0$ has the Levi decomposition:
$
P_0=(S\times G_{d_0})\ltimes N_0
$
where $N_0$ is the unipotent radical of $P_0$.
Also, with respect to the order of the basis in \eqref{bs}, the group $G_n$ is also defined by the following symmetric matrix:
\begin{equation} \label{eq:J}
J_{\Fr}^\Fn=\begin{pmatrix}
&&1\\&J_{\Fr-1}^{\Fn-2}&\\1&&
\end{pmatrix}_{\Fn\times\Fn}
\text{ and }
J_{0}^{\Fd_0}=\diag\{d_{1},\dots,d_{\Fd_0}\}
\end{equation}
as defined inductively.

\subsection{Discrete spectrum}\label{ssec-ds}
Following \cite{A13}, we denote by $\CA_2(G_n)$ the set of equivalence classes of irreducible
unitary representations $\pi$ of $G_n(\BA)$, where $\BA$ is the ring of adeles of the number field $k$, occurring in the discrete spectrum $L^2_\disc(G_n)$ of $L^2(G_n(F)\bks G_n(\BA))$. 
Also denote by $\CA_\cusp(G_n)$ for the subset of $\CA_2(G)$, whose elements occur in the cuspidal spectrum $L^2_\cusp(G_n)$.
The theory of endoscopic classification for $G_n$ is to parameterize the set $\CA_2(G_n)$ by means of the global Arthur parameters,
which can be realized as certain automorphic representations of general linear group $\GL_{2n}$. 

Let $G_n^\vee$ be the complex dual group of $G_n$, which is $\SO_{2n}(\BC)$ if $\Fn=2n$ and is $\Sp_{2n}$ if $\Fn=2n+1$. 
Following \cite{A13}, we denote by $\wt{\CE}_\simp(\GL_{2n})$ the set of the equivalence classes of simple twisted
endoscopic data of $\GL_{2n}$. Each member in $\wt{\CE}_\simp(\GL_{2n})$ is represented by a triple $(G_n^*,s,\xi)$, where $G_n^*$ is a $k$-quasisplit classical group,
$s$ is a semi-simple
element as described in \cite[Page 11]{A13}, and $\xi$ is the $L$-embedding: 
$
{^LG_n}\rightarrow {^L\GL_{2n}}.
$
The set of global Arthur parameters for $G_n^*$ is denoted by $\wt{\Psi}_2(G_n^*,\xi)$, or simply by $\wt{\Psi}_2(G_n^*)$ if the
$L$-embedding $\xi$ is well understood in the discussion. The elements of $\wt{\Psi}_2(G_n^*)$ can be formally written as 
\begin{equation}\label{aps}
\psi=(\tau_1,b_1)\boxplus\cdots\boxplus(\tau_r,b_r),
\end{equation}
where for each $j$, $\tau_j$ belongs to $\CA_\cusp(\GL_{a_j})$ and is self-dual, and $b_j$ is a positive integer. It follows that 
$2n=\sum_{j=1}^ra_jb_j$. A self-dual $\tau\in\CA_\cusp(\GL_a)$ is called {\it of symplectic type} if the (partial) exterior square $L$-function
$L^S(s,\tau,\wedge^2)$ has a (simple) pole at $s=1$; otherwise, $\tau$ is called {\it of orthogonal type}. In the latter case,
the (partial) symmetric square $L$-function $L^S(s,\tau,\sym^2)$ has a (simple) pole at $s=1$. 
A global parameter $\psi$ as in \eqref{aps} is called {\sl generic} if $b_i=1$ for $i=1,2,\cdots,r$.
The set of generic, elliptic, global Arthur parameters is denoted by $\wt{\Phi}_2(G_n^*)$. 
When $r=1$, the parameter $\psi=(\tau,b)$ is called {\sl simple}. 

\begin{thm}[Endoscopic Classification \cite{A13}]\label{ds}
For  each $\pi\in\CA_2(G_n)$, there is a $G_n$-relevant global Arthur parameter $\psi\in\wt{\Psi}_2(G_n^*,\xi)$,
such that $\pi$ belongs to the global Arthur packet, $\wt{\Pi}_{\psi}(G_n)$, attached to the global Arthur parameter $\psi$.
\end{thm}

Following \cite[Theorem 1.5.2]{A13},
when $G_n$ is an odd special orthogonal group, i.e. $\Fn=2n+1$, the multiplicity of $\pi\in\CA_2(G_n)$ realizing in the discrete spectrum $L_\disc^2(G_n)$ is
expected to be one. However, when $G_n$ is an even special orthogonal group, the discrete multiplicity of $\pi\in\CA_2(G_n)$ could be two.
In this case, we need to fix a realization of $\pi\in\CA_2(G_n)$ in the discrete spectrum $L_\disc^2(G_n)$, which will be denoted by $\CC_\pi$,
especially when the discrete multiplicity of $\pi$ is two.

Recall from \cite[Chapter 8]{A13} (see also \cite{JZ-BF}) that
\begin{equation}\label{wtO}
\wt{\RO}(G_n):=\wt{\out}_N(G_n):=\wt{\Aut}_N(G_n)/\wt{\Int}_N(G_n)
\end{equation}
is regarded as a diagonal subgroup of $\wt{\out}_N(G_n(\BA))$.
When $G_n$ is an even special orthogonal group,
one may take $\varepsilon\in \RO_{2n}(k)$ with $\det\varepsilon=-1$ and $\varepsilon^2=I_{2n}$, \label{pg:eps} such that
the action of $\wt{\RO}(G_n)$ on $\pi$ can be realized as the $\varepsilon$-conjugate on $\pi$,
i.e., $\pi^\varepsilon(g)=\pi(\varepsilon g\varepsilon^{-1})$. Hence the $\wt{\RO}(G_n)$-orbit of $\pi$ has one or two elements.
If $\wt{\RO}(G_n)$ acts freely on $\pi$, following the notation in \cite{A13}, we denote the $\wt{\RO}(G_n)$-torsor of $\pi$ by $\{\pi,\pi_\star\}$. \label{pg:star}
When $G_n$ is not an even special orthogonal group, the group $\wt{\RO}(G_n)$ is trivial, so is its action. Hence in this case,
the $\wt{\RO}(G_n)$-orbit of $\pi$ contains only $\pi$ itself.

When $G_n$ is an even special orthogonal group, an elliptic global Arthur parameter $\psi$ as in \eqref{aps} may be conjugated under the action of $\wt{\RO}(G_n)$ 
to a 
different global Arthur parameter $\psi_\star$ for $G_n$, which form an $\wt{\RO}(G_n)$-orbit $\{\psi,\psi_\star\}$. 
If the $\wt{\RO}(G_n)$-orbit is an $\wt{\RO}(G_n)$-torsor $\{\psi,\psi_\star\}$, then they define
different global Arthur packets and different global Vogan packets. In the rest of this paper, when we say that $\psi$ is a global Arthur parameter of an even special orthogonal group $G_n$,
we really mean that $\psi$ is identified with either $\psi$ or $\psi_\star$, through a specific twisted endoscopic datum.

\subsection{Bessel modules and Bessel Periods}\label{ssec-bmbp}
The general notion of Fourier coefficients of automorphic forms on $G_n(\BA)$ associated to nilpotent orbits was discussed in previous work (see \cite[Section 4]{J14} for instance). The Bessel-Fourier coefficient 
is a special type of Fourier coefficients associated to nilpotent that are parameterized by the following partition
\begin{equation}\label{bfp}
\udl{p}_\ell=[(2\ell+1)1^{\Fn-2\ell-1}].
\end{equation}
This partition is $G_n$-relevant if $0\leq\ell\leq\Fr$. In this case, 
we may choose the unipotent subgroup $V_\ell:=V_{\udl{p}_\ell}$ of $G_n$ to consist of all unipotent elements of the form:
\begin{equation}\label{vell}
V_{\ell}=\cpair{v=\begin{pmatrix}z&y&x\\&I_{\Fn-2\ell}&y'\\&&z^{*} \end{pmatrix}\in G_n \mid z\in Z_{\ell}},
\end{equation}
where $Z_{\ell}$ is the standard maximal (upper-triangular) unipotent subgroup of $\GL_\ell$. Note that
the $k$-rational nilpotent orbits $\CO_\ell$ in the $k$-stable nilpotent orbit $\CO^\st_{\udl{p}_\ell}$ are
in one-to-one correspondence with the $\GL_1\times G_{n-\ell}$-orbits of $k$-anisotropic vectors in
$(k^{\Fn-2\ell},q)$, which can be viewed as a subspace of $(V,q)$ and consists of the last row of $y$ in \eqref{nell}. It follows that the generic character $\psi_{\CO_\ell}$ of $V_\ell(\BA)$, 
whose restriction to $V_\ell(k)$, may also be explicitly defined as follows.
Fix a nontrivial character $\psi_k$ of $k\bks \BA$
and consider the following identification:
$$
V_{\ell}/[V_{\ell},V_{\ell}]\cong\oplus_{i=1}^{\ell-1}\Fg_{\alpha_{i}}\oplus k^{\Fn-2\ell}.
$$
Let $w_{0}$ be an anisotropic vector in $(k^{\Fn-2\ell},q)$ and define a character $\psi_{\ell,w_{0}}$ of $V_{\ell}(\BA)$ by
\begin{equation}\label{chw0}
\psi_{\CO_\ell}(v)=\psi_{\ell,w_{0}}(v):=\psi_k(\sum^{\ell-1}_{i=1}z_{i,i+1}+q(y_{\ell}, w_{0})),
\end{equation}
where $y_{\ell}$ is the last row of $y$ as defined in  \eqref{vell}.
We denote by $H_{\ell^-}^{\CO_\ell}=H_{\ell^-}^{w_0}$ is the identity connected component of the
stabilizer of  $\psi_{\ell,w_{0}}$ in $\GL_1\times G_{n-\ell}$, which is given by
\begin{equation}\label{Lellw0}
\cpair{\begin{pmatrix} I_{\ell}&&\\&\gamma&\\&&I_{\ell} \end{pmatrix}\in G_n \mid
\gamma J_{\Fn-2\ell}w_{0}=J_{\Fn-2\ell}w_{0}},
\end{equation}
where $\ell^-=[\frac{\Fl^-}{2}]$ with $\Fl^-:=\Fn-2\ell-1$. We may write $(k^{\Fn-2\ell},q)=(k^{\Fn-2\ell-1},q)\perp kw_0$. 
By \cite[Proposition 2.5]{JZ-BF}, the $k$-anisotropic kernel of the quadratic space $(k^{\Fn-2\ell-1},q)$ has dimension equal to $\Fd_0-1$ if $w_0$ belongs to the $\GL_1\times G_{n-\ell}$-orbit of a non-zero 
vector in the $k$-anisotropic kernel $(V_0,q)$ of $(V,q)$; otherwise, it has dimension equal to $\Fd_0+1$. This determines the structure of the group $H_{\ell^-}^{\CO_\ell}=H_{\ell^-}^{w_0}$.

Let $\varphi$ be an automorphic form on $G_n(\BA)$. We define the $\psi_{\ell,w_{0}}$-Bessel-Fourier coefficient of $\varphi$ by the following integral:
\begin{equation}\label{bf}
\CF^{\psi_{\ell,w_{0}}}(\varphi)(g):=
\int_{V_\ell(k)\bs V_\ell(\BA)}\varphi(vg)\psi_{\ell,w_{0}}^{-1}(v)dv.
\end{equation}
It is clear that $\CF^{\psi_{\ell,w_{0}}}(\varphi)(g)$ is a smooth automorphic function on $H_{\ell^-}^{w_0}$ with moderate growth. 
For any $\pi\in\CA_2(G_n)$, which has a realization $\CC_\pi$ in the discrete spectrum $L^2_\disc(G_n)$,
we denote by 
\begin{equation}\label{bm}
\CF^{\psi_{\ell,w_{0}}}(\CC_\pi)=\CF^{\CO_\ell}(\CC_\pi)
\end{equation}
the space generated by 
all $\CF^{\psi_{\ell,w_{0}}}(\varphi_\pi)(g)$ with all $\varphi_\pi\in \CC_\pi$,  
and call it an $\ell$-th {\sl Bessel module of $\pi$}. 

Let $\sigma\in\CA_2(H_m)$ have a realization $\CC_\sigma$ in the space of the discrete spectrum of $H_m(\BA)$. The Bessel period of $\CC_\pi$ and $\CC_\sigma$ is defined to be 
\begin{equation}\label{bp-L2}
\CB^{\CO_{\ell}}(\varphi_\pi,\varphi_\sigma)
:=
\int_{H_{\ell^-}^{\CO_{\ell}}(k)\bs H_{\ell^-}^{\CO_{\ell}}(\BA)}
\CF^{\CO_{\ell}}(\varphi_\pi)(h)\ol{\varphi_\sigma(h)}dh.
\end{equation}
Note that if both $\pi$ and $\sigma$ are not cuspidal, then the Bessel period $\CB^{\CO_{\ell}}(\varphi_\pi,\varphi_\sigma)$ needs regularization to define (see \cite{GJR04} \cite{GJR05} and \cite{GJR09} 
for instance). In spirit of the relation between parabolic induction and Bessel periods, we may only consider in this paper the situation that at least one of $\pi$, $\sigma$ is cuspidal. Hence the Bessel 
periods are well defined. 

In order to understand Question~\ref{q-L2bp}, the {\sl $L^2$-Automorphic Branching Problem}, and Question~\ref{q-L2rbp}, the {\sl $L^2$-Automorphic Reciprocal Branching Problem}, we would like to understand 
the Bessel periods $\CB^{\CO_{\ell}}(\varphi_\pi,\varphi_\sigma)$ in terms of other invariants attached to $\pi$ and $\sigma$. One natural invariant attached to $\pi$ and $\sigma$ are the relevant automorphic 
$L$-functions; and another natural invariant attached to $\pi$ and $\sigma$ is the global Arthur parameters. We will discuss the relation between the Bessel periods $\CB^{\CO_{\ell}}(\varphi_\pi,\varphi_\sigma)$ with 
the relevant automorphic $L$-functions in Section \ref{ssec-gzi}, and the relation between the Bessel periods $\CB^{\CO_{\ell}}(\varphi_\pi,\varphi_\sigma)$ with the relevant global Arthur parameters in 
Section \ref{sec-BPGAP}.

\subsection{A family of global zeta integrals}\label{ssec-gzi}
We recall from \cite[Section 4]{JZ-BF} a family of global zeta integrals that represent tensor product $L$-functions for the product of orthogonal groups with general linear groups. 
The global zeta integrals are nothing but the Bessel periods as defined in \eqref{bp-L2} with one of the cuspidal automorphic forms replaced by an Eisenstein series. 

Let $\tau$ be an irreducible unitary automorphic representation
of $\GL_a(\BA)$ of the following isobaric type:
\begin{equation}\label{tau8}
\tau=\tau_1\boxplus\tau_2\boxplus\cdots\boxplus\tau_r,
\end{equation}
where $\tau_i\in\CA_\cusp(\GL_{a_i})$, $\sum^{r}_{i=1}a_{i}=a$, and $\tau_i\not\cong\tau_j$ if $i\neq j$.
We may replace the last condition on $\tau_j$ by assuming that $\tau$ is generic, i.e. has a non-zero Whittaker-Fourier coefficient. 

Let $H_{a+m}$ be the special orthogonal group with a Levi subgroup $M_a=\GL_a\times H_m$. Since $(G_n,H_m)$ is a relevant pair in the sense of the Gan-Gross-Prasad conjecture (\cite{GGP12}), the 
pair $(G_n,H_{a+m})$ or $(H_{a+m},G_n)$ (depending on the size of $n$ and $a+m$) is also relevant. 
Let $P_a=M_a U_a$ be the standard parabolic $k$-subgroup of $H_{a+m}$ with the Levi subgroup $M_a=\GL_a\times H_m$. For $\sigma\in\CA_\cusp(H_m)$, we obtain the datum $(M_a,\tau\otimes\sigma)$, which is 
going to define an Eisenstein series, following \cite{L76} and \cite{MW95}. 
For any
\begin{equation}\label{af}
\phi=\phi_{\tau\otimes\sigma}\in \CA(U_a(\BA)M_a(F)\bks H_{a+m}(\BA))_{\tau\otimes\sigma},
\end{equation}
set $\phi_s:=|\det|^s_\BA\cdot\phi$, where $|\det|_\BA$ is first defined on $M_a(\BA)$, and then defined on $H_{a+m}(\BA)$ through the Iwasawa decomposition with respect to $P_a(\BA)$ and the standard maximal 
compact subgroup of $H_{a+m}(\BA)$. 
The associated Eisenstein series is defined to be
\begin{equation}\label{es}
E(h,\phi,s)=E(h,\phi_{\tau\otimes\sigma},s)=\sum_{\delta\in P_a(F)\bks H_{a+m}(F)}\phi_s(\delta g).
\end{equation}
The theory of Langlands on Eisenstein series (\cite{L76} and \cite{MW95}) shows that $E(h,\phi,s)$ converges absolutely for $\Re(s)$ large, has
meromorphic continuation to the complex plane $\BC$, and defines an automorphic form on $H_{a+m}(\BA)$ when $s$ is not a pole.

As in \eqref{bfp}, we take a family of $H_{a+m}$-relevant partitions 
$$
\udl{p}_\kappa=[(2\kappa+1) 1^{\Fm+2a-2\kappa-1}]
$$ 
for $H_{a+m}$ with
$\kappa\leq a+\Fr_\Fm$, where $\Fr_\Fm:=\Fr(H_m)$ is the $k$-rank of $H_m$. As in \eqref{bf}, we define the Bessel-Fourier coefficient of the Eisenstein series $E(h,\phi,s)$ on $H_{a+m}(\BA)$:
\begin{equation}\label{bfes}
\CF^{\psi_{\CO_\kappa}}(E(\cdot,\phi,s))(h)
:=
\int_{N_\kappa(k)\bks N_\kappa(\BA)}E(nh,\phi,s)\psi_{\CO_\kappa}^{-1}(n) \ud n,
\end{equation}
where the unipotent subgroup $N_\kappa$ of $H_{a+m}$ determined by the partition $\udl{p}_\kappa$ is similar to the unipotent subgroup
$V_\ell$ of $G_n$ as in \eqref{vell}, and the character $\psi_{\CO_\kappa}(n)$ is defined as in \eqref{chw0}. 
Hence its centralizer $G_{m^-}^{\CO_\kappa}$ is similar to the subgroup as described in \eqref{Lellw0}, with 
$m^-=[\frac{\Fm^-}{2}]$ and $\Fm^-:=2a+\Fm-2\kappa-1$ and $\kappa\leq a+\Fr_\Fm$. 
By \cite[Proposition 2.6]{JZ-BF},
the stabilizer $G_{m^-}^{\CO_\kappa}$ and the subgroup $H_m$ form a relevant pair in the sense of the
Gan-Gross-Prasad conjecture (\cite{GGP12}). In this case, the Bessel subgroup of $H_{a+m}$ is defined to be 
\begin{equation}\label{bsg}
R_{\CO_\kappa}^{a+m}:=G_{m^-}^{\CO_\kappa}\ltimes N_\kappa.
\end{equation}

For $\pi\in\CA_\cusp(G_{m^-}^{\CO_\kappa})$,
the {\sl global zeta integral} $\CZ(s,\phi_{\tau\otimes\sigma},\varphi_\pi,\psi_{\CO_\kappa})$, as in \cite[(4.9)]{JZ-BF},
is defined by the following Bessel period:
\begin{equation}\label{gzi-bp}
\CZ(s,\varphi_\pi,\phi_{\tau\otimes\sigma},\psi_{\CO_\kappa})
:=
\CB^{\CO_\kappa}(E(\phi_{\tau\otimes\sigma},s),\varphi_\pi).
\end{equation}
As given in Proposition~2.1 of \cite{JZ14}, $\CZ(s,\phi_{\tau\otimes\sigma},\varphi_\pi,\psi_{\CO_\kappa})$
converges absolutely and hence is holomorphic at $s$ where the Eisenstein series $E(h,\phi,s)$ has no poles. 
The main result  of \cite[Section 4]{JZ-BF} is the following theorem.

\begin{thm}[Basic Identity]\label{thm-bi}
With the notation given as above, the global zeta integral 
$\CZ(s,\phi_{\tau\otimes\sigma},\varphi_\pi,\psi_{\CO_\kappa})$ can be expressed as an eulerian product with factorizable data, and can be expressed by the following formula 
\begin{equation}\label{gzi-Lfn}
\CZ(s,\phi_{\tau\otimes\sig},\varphi_{\pi},\psi_{\CO_\kappa})
=\CZ_S(s,\cdot) \frac{L^S(s+\frac{1}{2},\tau\times\pi)}
{L^S(s+1,\tau\times\sig)L^S(2s+1,\tau,\rho)}
\end{equation}
where $\CZ_S(s,\cdot)$ is a finite eulerian product of the local zeta integrals at archimedean local places and ramified finite local places, $L^S(\cdot)$ denotes the partial $L$-function with unramified data, 
and 
$\rho=\wedge^2$ if $H_{m+a}$ is an even orthogonal group; and 
$\rho=\sym^2$ if $H_{m+a}$ is an odd orthogonal group. 
\end{thm}
Note that in the proof of the formula in \eqref{gzi-Lfn}, it is proved in \cite[Section 4]{JZ-BF} that the global zeta integral is not identically zero if and only if the Bessel period for 
$\pi$ and $\sig$ is not identically zero. By establishing some refined properties of the global intertwining operators and the local zeta integrals, we obtain the reciprocal relation of the relevant 
Bessel periods, which will be discussed in the next section.


\section{Bessel Periods and Global Arthur Parameters}\label{sec-BPGAP}


We are going to understand the relation between the Bessel periods for $\pi$ and $\sigma$ and the structure of the global Arthur parameters. The main tool is the basic identity \eqref{gzi-Lfn} in Theorem \ref{thm-bi}.
As explained in \cite[Section 5]{JZ-BF}, in order to find fruitful applications of the basic identity \eqref{gzi-Lfn}, one needs to have a good control of the analytic properties of the 
finite eulerian product of local zeta integrals, which is $\CZ_S(s,\cdot)$ on the right hand side of the basic identity; meanwhile, one needs to have a good control of the analytic properties of the Eisenstein series 
on the left hand side of the basic identity. Up to date, we are able to handle these two technical points when one of the $\sigma\in\CA_\cusp(H_m)$ and $\pi\in\CA_\cusp(G_n)$ 
has a generic global Arthur parameter. Hence we have to limit our discussion 
with such an {\sl assumption}, and hope to get back to those technical issues in our future work.

\subsection{$\pi$ has a generic global Arthur parameter}\label{ssec-pi-g}
Let $(G_n, H_m)$ be a relevant pair. 
Assume that  $\pi\in\CA_\cusp(G_n)$ has a generic global Arthur parameter. For any $\sigma\in\CA_\cusp(H_m)$, assume that for some $\ell$ with $0\leq\ell\leq\Fr_n=\Fr(G_n)$, 
the Bessel period $\CB^{\CO_{\ell}}(\varphi_\pi,\varphi_\sigma)$ is non-zero 
for a choice of $\varphi_\pi\in\CC_\pi$ and $\varphi_\sig\in\CC_\sig$. What can one say about the structure of the global Arthur parameter of $\sigma$?

By assumption, we write the generic global Arthur parameter of $\pi$ by 
\begin{equation}\label{gap-pi}
\phi_n=(\eta_1,1)\boxplus\cdots\boxplus(\eta_s,1). 
\end{equation}
We may write the global Arthur parameter of $\sigma$ by
\begin{equation}\label{ap-sigma}
\psi_m=(\zeta_1,2b_1+1)\boxplus\cdots\boxplus(\zeta_l,2b_l+1)\boxplus(\xi_1,2a_1)\boxplus\cdots\boxplus(\xi_k,2a_k).
\end{equation}
Note that $\eta_i$ and $\xi_j$ are in the same parity as irreducible cuspidal automorphic representations of general linear groups, i.e. either they are all of symplectic type or they are all of orthogonal type; 
and $\eta_i$ and $\zeta_j$ are in different parity, as irreducible cuspidal automorphic representations of general linear groups. Also the global Arthur parameter $\psi_m$ in \eqref{ap-sigma} includes the 
following three cases that 
$l=0$, but $k\neq 0$; $l\neq 0$, but $k=0$; and $l\times k\neq 0$. Note that $k=0$ is the same as $a_1=a_2=\cdots=a_k=0$.

 
\begin{thm}[Structure of $\psi_m$]  \label{thm-psim}
With $\pi$ and $\sigma$ as given above, and with the assumption that $\pi$ has a generic global Arthur parameter $\phi_n$ as in \eqref{gap-pi}, 
if the Bessel period $\CB^{\CO_{\ell}}(\varphi_\pi,\varphi_\sigma)$ is non-zero 
for a choice of $\varphi_\pi\in\CC_\pi$ and $\varphi_\sig\in\CC_\sig$, then the global Arthur parameter $\psi_m$ of $\sigma$ as given in \eqref{ap-sigma} must be of the form:
\begin{equation}\label{eq:psim}
\psi_m=(\zeta_1,1)\boxplus\cdots\boxplus(\zeta_l,1)\boxplus(\xi_1,2)\boxplus\cdots\boxplus(\xi_k,2),	
\end{equation}
and 
$\{\xi_1,\xi_2,\dots,\xi_k\}$ is a subset of $\{\eta_1,\eta_2,\dots,\eta_s\}$.
\end{thm}	
Note that the forms of $\phi_n$ in \eqref{gap-pi} and $\psi_m$ in \eqref{eq:psim} force that both $\{\xi_1,\xi_2,\dots,\xi_k\}$ and $\{\eta_1,\eta_2,\dots,\eta_s\}$ are sets instead of multiset. 
A special case when $\pi$ has a non-zero Whittaker-Fourier coefficiant, which implies that $G_n$ and $H_m$ must be quasi-split, was treated in \cite{S-1} and \cite{S-2} by local unramified argument. 

\begin{proof}
Without loss of generality, we may assume that 
\begin{equation}\label{a,b}
0\leq a_1\leq a_2\leq\cdots\leq a_k\ \  \textit{and}\ \ \  0\leq b_1\leq b_2\leq\cdots\leq b_l \text{ in \eqref{ap-sigma},}
\end{equation}
and if $k\neq 0$, then $a_1\geq 1$.
To prove this theorem, we consider two cases separately: {\it Case: $0\leq k$, $a_k\leq b_l$, and $l\neq 0$} and {\it Case: $a_k\geq b_l+1$; or $k\neq 0$ but $l=0$}.

{\bf Case: $0\leq k$, $a_k\leq b_l$, and $l\neq 0$.}\  In this case, we 
may assume that $\zeta_l\in\CA_\cusp(\GL_{\alpha_l})$. Note $\zeta_l$ is in different parity with $\eta_i$. Let $G_{\alpha_l+n}$ be the special orthogonal group containing a standard Levi subgroup 
$M_{\alpha_l}=\GL_{\alpha_l}\times G_n$. Then associated to the cuspidal datum $(M_{\alpha_l},\zeta_l\otimes\pi)$, we define an Eisenstein series $E(g,\phi_{\zeta_l\otimes\pi},s)$ as in \eqref{es}.
By Proposition 2.6 of \cite{JZ-BF}, we are able to choose an integer $\kappa$ with $0\leq\kappa\leq \alpha_l+\Fr_n=\Fr(G_n)$, and find a $k$-rational unipotent orbit $\CO_\kappa$ such that 
$\Fk^-=2\alpha_l+\Fn-2\kappa-1=\Fm$ with $\kappa^-=[\frac{\Fk^-}{2}]=m$, and 
the stabilizer $G_{\kappa^-}^{\CO_\kappa}=H_m$. As in \eqref{gzi-bp}, the Bessel period of $E(g,\phi_{\zeta_l\otimes\pi},s)$ and $\varphi_\sig$ gives a global zeta integral 
$\CZ(s,\phi_{\zeta_l\otimes\pi},\varphi_{\sig},\psi_{\CO_\kappa})$. By Theorem \ref{thm-bi}, we have the following identity:
\begin{equation}\label{gzi-Lfn2}
\CZ(s,\phi_{\zeta_l\otimes\pi},\varphi_{\sigma},\psi_{\CO_\kappa})
=\CZ_S(s,\cdot)\cdot\frac{L^S(s+\frac{1}{2},\zeta_l\times\sig)}
{L^S(s+1,\zeta_l\times\pi)L^S(2s+1,\zeta_l,\rho)}.
\end{equation}
Because the Bessel period  $\CB^{\CO_{\ell}}(\varphi_\pi,\varphi_\sigma)$ is non-zero 
for a choice of $\varphi_\pi\in\CC_\pi$ and $\varphi_\sig\in\CC_\sig$, the identity in \eqref{gzi-Lfn2} holds non-trivially. 

We first consider the right hand side of the basic identity \eqref{gzi-Lfn2}. By \eqref{ap-sigma}, the partial $L$-function $L^S(s+\frac{1}{2},\zeta_l\times\sig)$ can be written as 
\begin{eqnarray}
L^S(s+\frac{1}{2},\zeta_l\times\sig)
&=&\prod_{j=1}^lL^S(s+\frac{1}{2},\zeta_l\times (\zeta_j,2b_j+1))\nonumber\\
&&\cdot\prod_{i=1}^kL^S(s+\frac{1}{2},\zeta_l\times(\xi_i,2a_i)).
\end{eqnarray}
Then we look at each factors in the expression above, and have
\begin{equation}\label{zeta,zeta}
L^S(s+\frac{1}{2},\zeta_l\times (\zeta_j,2b_j+1))
=
\prod_{i=0}^{2b_j}L^S(s+\frac{1}{2}+b_j-i,\zeta_l\times\zeta_j),
\end{equation}
and
\begin{equation}\label{zeta,xi}
L^S(s+\frac{1}{2},\zeta_l\times(\xi_i,2a_i))
=
\prod_{j=0}^{2a_i-1}L^S(s+a_i-j,\zeta_l\times\xi_i).
\end{equation}
Since $\zeta_l$ is not equivalent to $\zeta_j$ for $j\neq l$, the right hand side of \eqref{zeta,zeta} has a ploe only when $j=l$ and the right-most (simple) pole is at $s=b_l+\frac{1}{2}$. 
If $k=0$, then the partial $L$-function $L^S(s+\frac{1}{2},\zeta_l\times\sig))$ has a right-most simple pole is at $s=b_l+\frac{1}{2}$.

If $k\neq 0$, then by assumption in this case, we have 
$$1\leq a_1\leq a_k\leq b_l.
$$
As irreducible cuspidal automorphic representations of general linear groups, $\zeta_j$ and $\xi_i$ are in different parity. It follows that the partial $L$-functions 
$L^S(s+a_i-j,\zeta_l\times\xi_i)$ has no poles. Since these $L$-functions are of symplectic type, they may have zero at $\Re(s)+a_i-j<1$, which may cancel the pole at $s=b_l+\frac{1}{2}$. 
We claim that this is impossible for all $i=1,2,\cdots,k$. 

In fact, if $\Re(s)+a_i-j<1$, then $\Re(s)<1+j-a_i$. Hence if the right-most simple pole at $s=b_l+\frac{1}{2}$ of $L^S(s+\frac{1}{2}-b_j,\zeta_l\times\zeta_l)$,
can be cancelled by the possible zero of $L^S(s+a_i-j,\zeta_l\times\xi_i)$ at  $\Re(s)+a_i-j<1$, we must have that 
$$
s=b_l+\frac{1}{2}<1+j-a_i.
$$
It follows that 
$$
a_i=1+(2a_i-1)-a_i\geq 1+j-a_i>b_l+\frac{1}{2}.
$$	
That is, $a_i> b_l+\frac{1}{2}$. 
Therefore we must have that  $a_k\geq b_l+1$. But this is impossible because the situation under consideration here assumes that $a_k\leq b_l$. 

This proves that the right-most simple pole 
at $s=b_l+\frac{1}{2}$ of the partial $L$-function in \eqref{zeta,zeta} can not be cancelled by any possible zeros of the partial $L$-functions in \eqref{zeta,xi}. 
Hence $L^S(s+\frac{1}{2},\zeta_l\times\sig)$ has a right-most simple pole at $s=b_l+\frac{1}{2}$.

By Proposition 5.5 of \cite{JZ-BF}, the finite eulerian product of local zeta integrals $\CZ_S(s,\cdot)$ at the right hand side of the basic identity \eqref{gzi-Lfn2} is a non-zero 
constant at $s=b_l+\frac{1}{2}$ for some choice of data, when the Bessel period  $\CB^{\CO_{\ell}}(\varphi_\pi,\varphi_\sigma)$ is non-zero. Therefore, the right hand side of \eqref{gzi-Lfn2} has a 
right-most simple at $s=b_l+\frac{1}{2}$ for some choice of data. 

Next we consider the left hand side of \eqref{gzi-Lfn}. The global zeta integral 
$\CZ(s,\phi_{\zeta_l\otimes\pi},\varphi_{\sig},\psi_{\CO_\kappa})$ is given by the Bessel period $\CB^{\CO_\kappa}(E(\cdot,\phi_{\zeta_l\otimes\pi},s),\varphi_\sig)$. 
According to the existence of the simple pole at $s=b_l+\frac{1}{2}$ on the right hand side, the Bessel period $\CB^{\CO_\kappa}(E(\cdot,\phi_{\zeta_l\otimes\pi},s),\varphi_\sig)$ has a simple pole at $s=b_l+\frac{1}{2}$. 
Since $\varphi_\sig$ is cuspidal, we obtain that the Eisenstein series $E(\cdot,\phi_{\zeta_l\otimes\pi},s)$ has a right-most simple pole at $s=b_l+\frac{1}{2}$. By Proposition 5.2 of \cite{JZ-BF}, 
if this Eisenstein series has a pole, then the pole must be at $s=\frac{1}{2}$. Hence we must have that $\frac{1}{2}=\frac{1}{2}+b_l$. This implies that $b_l=0$. Since the case under consideration assumes that 
$a_k\leq b_l$, we also have that $a_k=0$. Therefore, in this case, we must have that 
$$
a_1=a_2=\cdots=a_k=b_1=b_2=\cdots=b_l=0,
$$
which of course includes the case that $k=0$, and the global Arthur parameter of $\sigma$ by
$$
\psi_m=(\zeta_1,1)\boxplus\cdots\boxplus(\zeta_l,1).
$$

{\bf Case: $a_k\geq b_l+1$; or $k\neq 0$ but $l=0$.}\ In this case, we consider an Eisenstein series $E(\cdot,\phi_{\xi_k\otimes\pi},s)$. 
To set up the notation, we assume that $\xi_k$ is an irreducible cuspidal automorphic representation of $\GL_{\beta_k}(\BA)$. It follows 
that this Eisenstein series is defined on $G_{\beta_k+n}(\BA)$ with the defining cuspidal datum $(\GL_{\beta_k}\times G_n,\xi_k\otimes\pi)$. By Proposition 2.6 of \cite{JZ-BF}, there is a integer 
$\gamma$ with $0\leq\gamma\leq\beta_k+\Fr_n$, such that $H_m=G_{\Fq^-}^{\CO_\gamma}$, with $\Fq^-=2\beta_k+\Fn-2\gamma-1$. Then we consider the Bessel period 
$\CB^{\CO_\gamma}(E(\cdot,\phi_{\xi_k\otimes\pi},s),\varphi_\sig)$. As before, this Bessel period, which defines the global zeta integral for this case, produces the following basic identity: 
\begin{equation}\label{gzi-Lfn3}
\CB^{\CO_\gamma}(E(\cdot,\phi_{\xi_k\otimes\pi},s),\varphi_\sig)
=\CZ_S(s,\cdot)\cdot\frac{L^S(s+\frac{1}{2},\xi_k\times\sig)}
{L^S(s+1,\xi_k\times\pi)L^S(2s+1,\xi_k,\rho)}.
\end{equation}

In the right hand side of \eqref{gzi-Lfn3}, by \eqref{ap-sigma}, $L^S(s+\frac{1}{2},\xi_k\times\sig))$ can be written as 
\begin{eqnarray}
L^S(s+\frac{1}{2},\xi_k\times\sig)
&=&\prod_{j=1}^lL^S(s+\frac{1}{2},\xi_k\times (\zeta_j,2b_j+1))\nonumber\\
&&\cdot\prod_{i=1}^kL^S(s+\frac{1}{2},\xi_k\times(\xi_i,2a_i)). \label{eq:xik-sigma}
\end{eqnarray}
Each factors in the expression above can be further expressed as follows. 
\begin{equation}\label{xi,zeta}
L^S(s+\frac{1}{2},\xi_k\times (\zeta_j,2b_j+1))
=
\prod_{i=0}^{2b_j}L^S(s+\frac{1}{2}+b_j-i,\xi_k\times\zeta_j),
\end{equation}
and
\begin{equation}\label{xi,xi}
L^S(s+\frac{1}{2},\xi_k\times(\xi_i,2a_i))
=
\prod_{j=0}^{2a_i-1}L^S(s+a_i-j,\xi_k\times\xi_i).
\end{equation}
Since $\xi_k$ is not equivalent to $\xi_i$ for $i\neq k$, the right hand side of \eqref{xi,xi} has a pole only when $i=k$ and the right-most (simple) pole is at $s=a_k$. 
If $k\neq 0$ but $l=0$, the partial $L$-function $L^S(s+\frac{1}{2},\xi_k\times\sig)$ has a right-most (simple) pole is at $s=a_k$.

If $l\neq 0$, we have $a_k\geq b_l+1\geq 1$.
As irreducible cuspidal automorphic representations of general linear groups, $\zeta_j$ and $\xi_i$ are in different parity. It follows that the partial $L$-functions 
$L^S(s+a_i-j,\zeta_l\times\xi_i)$ has no poles. Since these $L$-functions are of symplectic type, they may have zero at $s+\frac{1}{2}+b_j-i=\frac{1}{2}$, which may cancel the pole at $s=a_k$. 
As before, we claim that this cancellation will not happen. 
If fact, if this cancellation happens, then we must have that $s=i-b_j=a_k$. Hence we have that $b_j\geq a_k$, in particular, we have that 
$b_l\geq a_k$. However, the situation under consideration here assumes that $a_k\geq b_l+1$. 

Therefore the partial $L$-function $L^S(s+\frac{1}{2},\xi_k\times\sig)$ must have a right-most simple pole at $s=a_k$, in the both cases considered in this part. 
Again, by Proposition 5.5 of \cite{JZ-BF}, the right hand side of \eqref{gzi-Lfn3} has a right-most simple pole at $s=a_k$, for some choice of data. It follows that in the left hand side, the Bessel period 
$\CB^{\CO_\gamma}(E(\cdot,\phi_{\xi_k\otimes\pi},s),\varphi_\sig)$ has a right-most simple pole at $s=a_k$ for some choice of data, and so does the Eisenstein series 
$E(\cdot,\phi_{\xi_k\otimes\pi},s)$. It is clear (\cite{JLZ13}) that the right-most possible pole of $E(\cdot,\phi_{\xi_k\otimes\pi},s)$ is at $s=1$. Hence, we must have that $a_k=1$, which forces that 
$b_l=0$. Therefore, we obtain that 
$$
a_1=a_2=\cdots=a_k=1\ \  \textit{and}\ \  b_1=b_2=\cdots=b_l=0,
$$
and 
the global Arthur parameter of $\sigma$ is of the form 
\begin{equation}\label{eq:psim-1}
\psi_m=(\zeta_1,1)\boxplus\cdots\boxplus(\zeta_l,1)\boxplus(\xi_1,2)\boxplus\cdots\boxplus(\xi_k,2).	
\end{equation}
Of course, if $l=0$, we must have that the global Arthur parameter of $\sigma$ is of the form 
$$
\psi_m=(\xi_1,2)\boxplus\cdots\boxplus(\xi_k,2).
$$

Now, let us assume that $\psi_m$ is of form \eqref{eq:psim-1}. 
We will show that 
$\{\xi_1,\xi_2,\dots,\xi_k\}$ is a subset of $\{\eta_1,\eta_2,\dots,\eta_s\}$.
For each $\xi_i$ for $1\leq i\leq k$, we consider the Bessel period
$\CB^{\CO_\gamma}(E(\cdot,\phi_{\xi_i}\otimes \pi, s),\varphi_\sig)$ and we have the same basic identity as \eqref{gzi-Lfn3}:
\begin{equation}\label{eq:gzi-Lfn3-xii}
\CB^{\CO_\gamma}(E(\cdot,\phi_{\xi_i\otimes\pi},s),\varphi_\sig)
=\CZ_S(s,\cdot)\cdot\frac{L^S(s+\frac{1}{2},\xi_i\times\sig)}
{L^S(s+1,\xi_i\times\pi)L^S(2s+1,\xi_i,\rho)}.
\end{equation}
Similar as above, $\CZ_S(s,\cdot)$ at the right hand side of \eqref{eq:gzi-Lfn3-xii} is nonzero at $s=1$ for some choice of data. 
Following \eqref{eq:xik-sigma}, \eqref{xi,zeta} and \eqref{xi,xi}, 
$L^S(s+\frac{1}{2},\xi_i\times\sig)$ has a simple pole at $s=1$.
Then the right hand side of \eqref{eq:gzi-Lfn3-xii} has a right-most simple pole at  $s=1$.

It follows that in the left hand side, the Bessel period
$\CB^{\CO_\gamma}(E(\cdot,\phi_{\xi_i}\otimes \pi, s),\varphi_\sig)$ has a right-most simple pole at  $s=1$ for some choice of data, 
so does the Eisenstein series $E(\cdot,\phi_{\xi_i}\otimes \pi, s)$.
Recall that $E(\cdot,\phi_{\xi_i}\otimes \pi, s)$ has a simple pole at $s=1$ if and only if $L^S(s,\xi_i\times \pi)$ has a simple pole at $s=1$.
Thus the existence of the simple pole of $L^S(s,\xi_i\times \pi)$ at $s=1$ 
is equivalent to  $\xi_i=\eta_j$ for some $\eta_j$.
Hence, we obtain that $\{\xi_1,\xi_2,\dots,\xi_k\}$ is a subset of $\{\eta_1,\eta_2,\dots,\eta_s\}$.

This finishes the proof. 
\end{proof}

Based on the above discussion, we make the following conjecture. 
\begin{conj}\label{conj-pi-generic}
Assume that $\pi\in\CA_\cusp(G_n)$ has a generic global Arthur parameter $\phi_n$ as in \eqref{gap-pi}. For any $\sigma\in\CA_\cusp(H_m)$, assume that the Bessel period 
$\CB^{\CO_{\ell}}(\varphi_\pi,\varphi_\sigma)$ is non-zero 
for a choice of $\varphi_\pi\in\CC_\pi$ and $\varphi_\sig\in\CC_\sig$. If the global Arthur parameter $\psi_m$ of the $\sigma$ has the part $(\xi_1,2)\boxplus\cdots\boxplus(\xi_k,2)$ with non-zero $k$, 
then the integer $\ell$ with $0\leq\ell\leq\Fr_n=\Fr(G_n)$ is not the first occurrence index of $\pi$. 
\end{conj}

\subsection{$\sigma$ has a generic global Arthur parameter}\label{ssec-sigma-g}
Let $G_n, H_m$ be a relevant pair. 
Assume that  $\sigma\in\CA_\cusp(H_m)$ has a generic global Arthur parameter. For any $\pi\in\CA_\cusp(G_n)$, assume that for some $\ell$ with $0\leq\ell\leq\Fr_n=\Fr(G_n)$, 
the Bessel period $\CB^{\CO_{\ell}}(\varphi_\pi,\varphi_\sigma)$ is non-zero 
for a choice of $\varphi_\pi\in\CC_\pi$ and $\varphi_\sig\in\CC_\sig$. What can one say about the structure of the global Arthur parameter of $\pi$? 
We may regard this question reciprocal to that in Section \ref{ssec-pi-g}. Indeed, we are going to apply the same theory developed in \cite[Sections 4,5]{JZ-BF}, but to different situation. 
For completeness, we include the details for this case. 

In this case, we may write the global Arthur parameters of $\pi$ and $\sigma$ as follows:
\begin{equation}\label{ap-pi}
\psi_n=(\eta_1,2a_1+1)\boxplus\cdots\boxplus(\eta_k,2a_k+1)\boxplus(\delta_1,2b_1)\boxplus\cdots\boxplus(\delta_l,2b_l)
\end{equation}
for $\pi$, and 
\begin{equation}\label{gap-sigma}
\phi_m=(\zeta_1,1)\boxplus\cdots\boxplus(\zeta_s,1)
\end{equation}
for $\sigma$. Note that $\delta_i$ and $\zeta_j$ are in the same parity, while $\eta_i$ and $\zeta_j$ are in different parity, as irreducible cuspidal automorphic representations of general linear groups. 
Also the global Arthur parameter $\psi_n$ in \eqref{ap-pi} includes the 
following three cases that 
$l=0$, but $k\neq 0$; $l\neq 0$, but $k=0$; and $l\times k\neq 0$. Note that $l=0$ is the same as $b_1=\cdots=b_l=0$. 

{\bf Cases: $0\leq l$, $b_l\leq a_k$ and $k\neq 0$.}\ 
We consider the Eisenstein series $E(\cdot,\phi_{\eta_k\otimes\sigma},s)$ of $H_{\alpha_k+m}(\BA)$, where we may assume that $\eta_k\in\CA_\cusp(\GL_{\alpha_k})$. It follows that the Eisenstein series 
has its cuspidal support $(M_{\alpha_k},\eta_k\otimes\sigma)$, with $M_{\alpha_k}=\GL_{\alpha_k}\times H_m$ being a standard Levi subgroup of $H_{\alpha_k+m}$. The global zeta integral that is defined to 
be the Bessel period $\CB^{\CO_\kappa}(E(\cdot,\phi_{\eta_k\otimes\sigma},s),\varphi_\pi)$ has the following identity, as proved in Section 4 of \cite{JZ-BF}, 
\begin{equation}\label{gzi-Lfn4}
\CB^{\CO_\kappa}(E(\cdot,\phi_{\eta_k\otimes\sig},s),\varphi_\pi)
=\CZ_S(s,\cdot)\cdot\frac{L^S(s+\frac{1}{2},\eta_k\times\pi)}
{L^S(s+1,\eta_k\times\sig)L^S(2s+1,\eta_k,\rho)}.
\end{equation}
We consider the partial $L$-function $L^S(s+\frac{1}{2},\eta_k\times\pi)$ on the right hand side: 
\begin{eqnarray}
L^S(s+\frac{1}{2},\eta_k\times\pi)
&=&\prod_{j=1}^kL^S(s+\frac{1}{2},\eta_k\times (\eta_j,2a_j+1))\nonumber\\
&&\cdot\prod_{i=1}^lL^S(s+\frac{1}{2},\eta_k\times(\delta_i,2b_i)).
\end{eqnarray}
As before, we write
\begin{equation}\label{eta,eta}
L^S(s+\frac{1}{2},\eta_k\times (\eta_j,2a_j+1))
=
\prod_{i=0}^{2a_j}L^S(s+\frac{1}{2}+a_j-i,\eta_k\times\eta_j),
\end{equation}
and
\begin{equation}\label{eta,delta}
L^S(s+\frac{1}{2},\eta_k\times(\delta_i,2b_i))
=
\prod_{j=0}^{2b_i-1}L^S(s+b_i-j,\eta_k\times\delta_i).
\end{equation}
Hence we obtain that the right-most simple pole of the partial $L$-function in \eqref{eta,eta} is at $s=a_k+\frac{1}{2}$. Hence if $l=0$, then the partial $L$-function 
$L^S(s+\frac{1}{2},\eta_k\times\pi)$ has a right-most simple pole at $s=a_k+\frac{1}{2}$.

If $l\neq 0$, since $\eta_j$ and $\delta_i$ are in different parity, the partial $L$-functions $L^S(s+b_i-j,\eta_k\times\delta_i)$ are of symplectic type. Hence their possible zero at 
$\Re(s)+b_i-j<1$, i.e. $\Re(s)<1+j-b_i$ may cancel the simple pole at $s=a_k+\frac{1}{2}$. If this is true, then we must have that $b_i> \frac{1}{2}+a_k$, in particular, that 
$b_l\geq 1+a_k$. But this contradicts the assumption that $b_l\leq a_k$ in this case. 

This simply proves that the partial $L$-function $L^S(s+\frac{1}{2},\eta_k\times\pi)$ has a right-most simple pole at $s=a_k+\frac{1}{2}$. By Proposition 5.5 of \cite{JZ-BF}, we obtain that 
the right hand side of the identity in \eqref{gzi-Lfn4} has a right-most simple pole at  $s=a_k+\frac{1}{2}$ for some choice of data. It follows that the left hand side of the identity in \eqref{gzi-Lfn4} has 
a right-most simple pole at $s=a_k+\frac{1}{2}$ for some choice of data, and so is the Eisenstein series $E(\cdot,\phi_{\eta_k\otimes\sig},s)$. Finally by Proposition 5.2 of \cite{JZ-BF}, we obtain that 
the Eisenstein series can have the right-most simple pole at $s=\frac{1}{2}$. Therefore, we obtain that $a_k=0$, which implies that $b_l=0$. This of course implies that 
$$
a_1=\cdots=a_k=b_1=\cdots=b_l=0.
$$
In particular, the global Arthur parameter $\psi_n$ of $\pi$ must be of the form:
\begin{equation}\label{psin-1}
\psi_n=(\eta_1,1)\boxplus\cdots\boxplus(\eta_k,1)
\end{equation}

{\bf Case: $b_l\geq a_k+1$; or $l\neq 0$ but $k=0$.}\ 
In this case, we consider the Eisenstein series $E(\cdot, \phi_{\delta_l\otimes\sigma},s)$, where we may assume that $\delta_l\in\CA_\cusp(\GL_{\beta_l})$. 
This Eisenstein series is defined over $H_{\beta_l+m}(\BA)$, with its cuspidal support $(M_{\beta_l},\delta_l\otimes\sigma)$, where $M_{\beta_l}=\GL_{\beta_l}\times H_m$ is a standard Levi subgroup of 
$H_{\beta_l+m}$. The Bessel period $\CB^{\CO_\gamma}(E(\cdot, \phi_{\delta_l\otimes\sigma},s), \varphi_\pi)$ defines the global zeta integral and has the following identity:
\begin{equation}\label{gzi-Lfn5}
\CB^{\CO_\gamma}(E(\cdot,\phi_{\delta_l\otimes\sig},s),\varphi_\pi)
=\CZ_S(s,\cdot)\cdot\frac{L^S(s+\frac{1}{2},\delta_l\times\pi)}
{L^S(s+1,\delta_l\times\sig)L^S(2s+1,\delta_l,\rho)},
\end{equation}
as proved in Section 4 of \cite{JZ-BF}. Following the same idea, we consider the partial $L$-functions on the right hand side of the identity in \eqref{gzi-Lfn5} and have 
\begin{eqnarray}
L^S(s+\frac{1}{2},\delta_l\times\pi)
&=&\prod_{j=1}^kL^S(s+\frac{1}{2},\delta_l\times (\eta_j,2a_j+1))\nonumber\\
&&\cdot\prod_{i=1}^lL^S(s+\frac{1}{2},\delta_l\times(\delta_i,2b_i)).
\end{eqnarray}
As before, we write
\begin{equation}\label{delta,eta}
L^S(s+\frac{1}{2},\delta_l\times (\eta_j,2a_j+1))
=
\prod_{i=0}^{2a_j}L^S(s+\frac{1}{2}+a_j-i,\delta_l\times\eta_j),
\end{equation}
and
\begin{equation}\label{delta,delta}
L^S(s+\frac{1}{2},\delta_l\times(\delta_i,2b_i))
=
\prod_{j=0}^{2b_i-1}L^S(s+b_i-j,\delta_l\times\delta_i).
\end{equation}
It is clear that the right-most simple pole of \eqref{delta,delta} is at $s=b_l$. If $k=0$, then the partial $L$-function $L^S(s+\frac{1}{2},\delta_l\times\pi)$ has a right-most simple pole at 
$s=b_l$. 

If $k\neq 0$, the partial $L$-functions $L^S(s+\frac{1}{2}+a_j-i,\delta_l\times\eta_j)$ are of symplectic type. Hence they might have a zero at $\Re(s)+\frac{1}{2}+a_j-i<1$. If this possible zero is going to 
cancel the right-most simple pole at $s=b_l$, we must have that $b_l+\frac{1}{2}+a_j-i<1$, which implies that $b_l-a_k<\frac{1}{2}$. But this is impossible since in this case, we assume that $b_l\geq a_k+1$. 

This proves that the partial $L$-function $L^S(s+\frac{1}{2},\delta_l\times\pi)$ has a right-most simple pole at $s=b_l$. By Proposition 5.5 of \cite{JZ-BF}, 
the right hand side of the identity in \eqref{gzi-Lfn5} has a right-most simple pole at  $s=b_l$ for some choice of data. Now, we apply this to the left hand side and obtain that the Bessel period 
$\CB^{\CO_\gamma}(E(\cdot, \phi_{\delta_l\otimes\sigma},s), \varphi_\pi)$ has a right-most simple pole at  $s=b_l$ for some choice of data, and so does the Eisenstein series 
$E(\cdot, \phi_{\delta_l\otimes\sigma},s)$. By \cite{JLZ13}, this Eisenstein seris can have a right-most simple pole at $s=1$. It follows that $b_l=1$ and $a_k=0$. 
It follows that if $k=0$, then the global Arthur parameter of $\pi$ is of the form:
$$
\psi_n=(\delta_1,2)\boxplus\cdots\boxplus(\delta_l,2);
$$
and if $k\neq 0$, then the global Arthur parameter of $\pi$ is of the form:
\begin{equation}\label{eq:psin}
\psi_n=(\eta_1,1)\boxplus\cdots\boxplus(\eta_k,1)\boxplus(\delta_1,2)\boxplus\cdots\boxplus(\delta_l,2).	
\end{equation}

Let us assume that $\psi_n$ is of form \eqref{eq:psin}.
Analogous to Theorem \ref{thm-psim}, we will show that  $\{\delta_1,\delta_2,\dots,\delta_l\}$ is a subset of $\{\zeta_1,\zeta_2,\dots,\zeta_s\}$.

Similar to the proof of Theorem \ref{thm-psim}
for {\it Case: $a_k\geq b_l+1$; or $k\neq 0$ but $l=0$},
we may consider the Bessel period $\CB^{\CO_\gamma}(E(\cdot, \phi_{\delta_i\otimes\sigma},s), \varphi_\pi)$ for $1\leq i\leq l$.
By computing $L^S(s+\frac{1}{2}, \delta_i\times \pi)$, 
we have the Bessel period has a right-most simple pole at $s=1$.
It follows that the Eisenstein series $E(\cdot, \phi_{\delta_i\otimes\sigma},s)$ has a  simple pole at $s=1$,
which implies that $L^S(s,\delta_i\times\pi)$ has a simple pole at $s=1$.
Hence $\delta_i=\zeta_j$ for some $1\leq j\leq s$.

We summarize the above discussion as a theorem.

\begin{thm}[Structure of $\psi_n$]\label{thm-psin}
With $\pi$ and $\sigma$ as given above, and with the assumption that $\sigma$ has a generic global Arthur parameter $\phi_m$ as in \eqref{gap-sigma}, 
if the Bessel period $\CB^{\CO_{\ell}}(\varphi_\pi,\varphi_\sigma)$ is non-zero 
for a choice of $\varphi_\pi\in\CC_\pi$ and $\varphi_\sig\in\CC_\sig$, then the global Arthur parameter $\psi_n$ of $\pi$ as given in \eqref{ap-pi} must be of the form:
$$
\psi_n=(\eta_1,1)\boxplus\cdots\boxplus(\eta_k,1)\boxplus(\delta_1,2)\boxplus\cdots\boxplus(\delta_l,2)
$$
and  $\{\delta_1,\delta_2,\dots,\delta_l\}$ is a subset of $\{\zeta_1,\zeta_2,\dots,\zeta_s\}$.
\end{thm}

A special case when $\sigma$ has a non-zero Whittaker-Fourier coefficients, which implies that $G_n$ and $H_m$ must be quasi-split, was treated in \cite{S-1} and
\cite{S-2} by local unramified argument. 


\subsection{Reformulation in terms of $L$-functions}\label{ssec-rlfn}
We continue our discussion for the pair $(\pi,\sigma)$ with one of them having a generic global Arthur parameter. We may first consider the situation that $\pi$ has a generic global Arthur parameter 
$\phi_n=(\eta_1,1)\boxplus\cdots\boxplus(\eta_s,1)$ as in \eqref{gap-pi}, and $\sigma$ has a global Arthur parameter $\psi_m$ as in \eqref{ap-sigma}. 
In the spirit of the global Gan-Gross-Prasad conjeccture 
for non-generic global Arthur parameters, as explained in \cite{GGP18}, we prove the following result.


 
\begin{thm}[Bessel Period and Central $L$-value]\label{thm-bpclv}
Assume that $\pi$ and $\sigma$ are given as given in Theorem \ref{thm-psim}, with $\pi$ having a generic global Arthur parameter $\phi_n$ as in \eqref{gap-pi} and $\sigma$ having a global Arthur parameter $\psi_m$ as 
\eqref{ap-sigma}, which is of the form: 
$$
\psi_m=(\zeta_1,2b_1+1)\boxplus\cdots\boxplus(\zeta_l,2b_l+1)\boxplus(\xi_1,2a_1)\boxplus\cdots\boxplus(\xi_k,2a_k).
$$
Assume that the following non-vanishing 
\begin{equation}\label{eq:nonzero-condition}
\prod_{1\leq i\leq l,~ 1\leq j\leq k}L(\frac{1}{2},\zeta_i\times\xi_j)\ne 0,
\end{equation}
holds. 
If the Bessel period $\CB^{\CO_{\ell}}(\varphi_\pi,\varphi_\sigma)$ is non-zero 
for a choice of $\varphi_\pi\in\CC_\pi$ and $\varphi_\sig\in\CC_\sig$
then 
\begin{equation}\label{eq:adjoint}
\frac{L(s+\frac{1}{2},\phi_{n}\times \psi_m)}{L(s+1,\phi_n,\Ad_G)L(s+1,\psi_m,\Ad_H)}
\end{equation}
is nonzero at $s=0$.
\end{thm}	

Note that here we use the complete $L$-functions of the Arthur parameters, which are defined in terms of the relevant complete $L$-functions associated to irreducible cuspidal automorphic representations of 
general linear groups. 

\begin{proof}
Since the assumption of the theorem is taken from Theorem \ref{thm-psim}, it is clear by Theorem \ref{thm-psim} that the global Arthur parameters $\phi_n$ and $\psi_m$ are of the form: 
\[
\phi_n=(\eta_1,1)\boxplus\cdots\boxplus(\eta_s,1),
\]
and
\[
\psi_m=(\zeta_1,1)\boxplus\cdots\boxplus(\zeta_l,1)\boxplus(\xi_1,2)\boxplus\cdots\boxplus(\xi_k,2)
\]
with the property that $\{\xi_1,\xi_2,\dots,\xi_k\}$ is a subset of $\{\eta_1,\eta_2,\dots,\eta_s\}$. We may 
denote by $\{\xi'_1,\xi'_2,\dots,\xi'_{k'}\}$ its complementary set. The global Arthur parameter $\phi_n$ yields an endoscopy structure $\CG_1\times\CG_2\times\cdots\times \CG_s$ of $G_n$
and $\psi_m$ yields an endoscopy structure $\CH_1\times\cdots\times\CH_l\times\CF_1\times\cdots\CF_k$ of $H_m$.

Since $\eta_i$, $\xi_j$ and $\zeta_{r}$ are self-dual, we may write 
\begin{equation}\label{phin-ad}
L(s,\phi_n,\Ad_{G_n})=\prod_{i=1}^{s}L(s,\eta_i,\Ad_{\CG_i})\prod_{1\leq i<j\leq s}L(s,\eta_i\times\eta_j),
\end{equation}
\begin{eqnarray}\label{psim-ad}
L(s,\psi_m,\Ad_{H_m})&=&\prod_{i=1}^{l}L(s,\zeta_i,\Ad_{\CH_i})\prod_{\substack{1\leq i \leq l\\ 1\leq j\leq k}}L(s+\frac{1}{2},\zeta_i\times\xi_j)L(s-\frac{1}{2},\zeta_i\times\xi_j)\nonumber\\
&&\prod_{j=1}^{k}L(s,\xi_j\times\xi_j)L(s+1,\xi_j,\rho)L(s-1,\xi_j,\rho) 
\end{eqnarray}
and 
\begin{equation}\label{phin-psim}
L(s,\phi_n\times \psi_m)=\prod_{\substack{1\leq i\leq l \\ 1\leq j \leq s}}
L(s, \zeta_i\times \eta_j)
\prod_{\substack{1\leq i \leq s\\ 1\leq j\leq k}}
L(s+\frac{1}{2},\eta_i\times\xi_j)L(s-\frac{1}{2},\eta_i\times\xi_j). 
\end{equation}
We are going to show that 
$$
\frac{L(s+\frac{1}{2},\phi_{n}\times \psi_m)}{L(s+1,\phi_n,\Ad_{G_n})L(s+1,\psi_m,\Ad_{H_m})}
$$
is non-zero at $s=0$. 

By \eqref{phin-ad}, $L(s+1,\phi_n,\Ad_{G_n})$ is holomorphic and non-zero at $s=0$, since $\eta_i$ are cuspidal and distinct to each other. It reduces to show that 
\begin{equation}\label{nonzero-1}
\frac{L(s+\frac{1}{2},\phi_{n}\times \psi_m)}{L(s+1,\psi_m,\Ad_{H_m})}
\end{equation}
is non-zero at $s=0$. 

To consider $L(s+1,\psi_m,\Ad_{H_m})$ at $s=0$, we study the equation in \eqref{psim-ad}. 
It is clear that at $s=0$, $L(s+1,\zeta_i,\Ad_{\CH_i})$, $L(s+\frac{3}{2},\zeta_i\times\xi_j)$, $L(s+\frac{1}{2},\zeta_i\times\xi_j)$, and $L(s+2,\xi_j,\rho)$ are holomophic and non-zero at $s=0$. 
The non-vanishng of \eqref{nonzero-1} at $s=0$ is equivalent to the non-vanishing of 
\begin{equation}\label{nonzero-2}
\frac{\prod_{\substack{1\leq i\leq l \\ 1\leq j \leq s}}L(s+\frac{1}{2}, \zeta_i\times \eta_j)
\prod_{\substack{1\leq i \leq s\\ 1\leq j\leq k}}L(s+1,\eta_i\times\xi_j)L(s,\eta_i\times\xi_j)}
{\prod_{\substack{1\leq i \leq l\\ 1\leq j\leq k}}L(s+\frac{1}{2},\zeta_i\times\xi_j)\prod_{j=1}^{k}L(s+1,\xi_j\times\xi_j)L(s,\xi_j,\rho)}
\end{equation}
at $s=0$. It is clear that 
$$
\frac{\prod_{\substack{1\leq i \leq s\\ 1\leq j\leq k}}L(s+1,\eta_i\times\xi_j)L(s,\eta_i\times\xi_j)}
{\prod_{j=1}^{k}L(s+1,\xi_j\times\xi_j)L(s,\xi_j,\rho)}
$$
is non-zero at $s=0$. It is sufficient to show that 
\begin{equation}\label{nonzero-3}
\prod_{1\leq i\leq l,~ 1\leq j \leq k'}
L(s+\frac{1}{2}, \zeta_i\times \xi'_j)		
=
\frac{\prod_{\substack{1\leq i\leq l \\ 1\leq j \leq s}}L(s+\frac{1}{2}, \zeta_i\times \eta_j)}
{\prod_{\substack{1\leq i \leq l\\ 1\leq j\leq k}}L(s+\frac{1}{2},\zeta_i\times\xi_j)}
\end{equation}
is non-zero at $s=0$.

Under the assumption \eqref{eq:nonzero-condition}, \eqref{nonzero-3} is equivalent to 
\[
\prod_{1\leq i\leq l,~ 1\leq j \leq s}
L(\frac{1}{2}, \zeta_i\times \eta_j)\ne 0.		
\] 
We are going to show that for each $\zeta_i$ ($1\leq i\leq l$) 
$$L(\frac{1}{2},\zeta_i\times\pi)=\prod_{1\leq j\leq s}L(\frac{1}{2}, \zeta_i\times \eta_j )$$ is nonzero. 

Consider the basic identity as \eqref{gzi-Lfn3} for the Bessel period 
$$\CB^{\CO_\gamma}(E(\cdot,\phi_{\zeta_i\otimes \pi}, s),\varphi_\sig).$$ 
The right hand side contains a partial $L$-function:
\begin{align*}
 &L^S(s+\frac{1}{2},\zeta_i\times\sig)\\
=&\prod_{j=1}^{l}L^S(s+\frac{1}{2},\zeta_i\times\zeta_j)\prod_{r=1}^{k}
L^S(s+1,\zeta_i\times\xi_r)L^S(s,\zeta_i\times\xi_r).
\end{align*}
Since $L_S(\frac{1}{2},\zeta_i\times\xi_r)=\prod_{\nu\in S}L(\frac{1}{2},\zeta_{i,\nu}\times\xi_{r,\nu})$ for $1\leq r\leq k$ is nonzero,
we have $L^S(\frac{1}{2},\zeta_i\times\xi_r)\ne 0$ by the assumption \eqref{eq:nonzero-condition}.
It follows that $L^S(s+\frac{1}{2},\zeta_i\times\sig)$ has a simple pole at $s=\frac{1}{2}$. 
Hence the right hand side of \eqref{gzi-Lfn3} for $\zeta_i$ has a simple pole at $s=\frac{1}{2}$ for some choice of date, so does the left hand side of \eqref{gzi-Lfn3} for $\zeta_i$.
It implies that the Eisenstein series $E(\cdot,\phi_{\zeta_i\otimes \pi}, s)$
has a simple pole at $s=\frac{1}{2}$. Since $\pi$ has ageneric global Arthur parameter $\phi_n$, by Proposition 5.2 of \cite{JZ-BF}, 
$$L(2s,\zeta_i,\rho)L(s,\zeta_i\times\pi)$$ 
must have  a simple pole at $s=\frac{1}{2}$.
Hence we must have that 
$$\prod_{1\leq j\leq s}L(s,\zeta_i\times \eta_j)\ne 0$$
for each $\zeta_i$. This completes the proof. 

\end{proof}

It is clear that the same result holds when $\sigma$ has a generic global Arthur parameter $\phi_m$ and $\pi$ has a global Arthur parameter $\psi_n$, by the same argument. We omit the statement here. 
Due to the technical issues related to the apllication of the basic identity (Theorem \ref{thm-bi}), we will not go on to discuss the more general situation and leave it for our future work.


\section{Theory of Twisted Automorphic Descents}\label{sec-TTAD}


We are going to explain the theory of twisted automorphic descents that mainly developed by the authors in \cite{JZ-BF}, as part of the general theory of explicit constructions of integral transforms 
with automorphic kernel functions to realize the endoscopic transfers and endoscopic descents as outlined in \cite{J14}. This theory is part of the effort to develop refined structures of the discrete spectrum 
of classical groups based the theory of endoscopic classification that has been achieved through the Arthur-Selberg trace formula approach. 

Although the idea of construction of automorphic descents goes back to 
\cite{An79}, \cite{PS83} and \cite{EZ85}, or maybe even earlier, the recent development of the automorphic descents was re-initiated by the work of Ginzburg-Rallis-Soudry, which has been well presented in their 
book \cite{GRS11}. The motivation of the work of Ginzburg-Rallis-Soudry is to construct the backwards Langlands functorial transfer from the general linear group $\GL_N$ to a quasi-split classical group $G$, 
corresponding to the standard embedding of the complex dual groups: $G^\vee(\BC)\rightarrow\GL_N(\BC)$, where $N$ depends on $G$.  
It is well-known that the existence of such Langlands functorial transfers can be proved either through 
the Arthur-Selberg trace formula approach, or through the approach that is based on the converse theorem of Cogdell and Piatetski-Shapiro. 

However, in order to develop a refined theory of discrete spectrum in general, it may be more important to understand the automorphic descent of Ginzburg-Rallis-Soudry as an explicit construction of 
cuspidal automorphic representations for classical groups. The explicit cuspidal modules on quasi-split classical groups that have 
been constructed through the Ginzburg-Rallis-Soudry descent have a non-zero Whittaker-Fourier coefficient. Hence it remains very interesting to ask: 
\begin{ques}\label{q-non-generic}
How to construct explicit cuspidal modules for cuspidal automorphic representations that has no non-zero Whittaker-Fourier coefficients? 
\end{ques} 
and even more generally, to ask: 
\begin{ques}\label{q-nonquasisplit}
How to construct cuspidal modules for the whole cuspidal spectrum when the classical groups are not quasi-split?
\end{ques}
In order to understand these questions, it is better to understand the construction theory in terms of the endoscopic classification of the discrete spectrum of classical groups. 
By Theorem \ref{ds}, any $\pi\in\CA_2(G)$ belong to a global Arthur packet $\wt{\Pi}_\psi(G)$ for some global Arthur parameter $\psi\in\wt{\Psi}_2(G)$. It is reasonable to ask 
\begin{ques}\label{q-basedonec}
How to construct an explicit cuspidal module for $\pi$ if assume that $\pi\in\wt{\Pi}_\psi(G)$?
\end{ques}
These questions have been discussed in \cite{JZ-BF}, related to the Principle of the Construction of cuspidal modules (\cite[Principle 1.1]{JZ-BF}). Instead of getting into the details here, we would like to recall 
what we have established in the theory of twisted automorphic descents. 

For more precise discussions, we only consider special orthogonal group $G_n$. 
The Ginzburg-Rallis-Soudry descent is to use the source automorphic representations of quasi-split $G_n^*$ with global Arthur parameters $\phi_\tau^{(2)}$, which is given by 
$$
\phi_\tau^{(2)}=(\tau_1,2)\boxplus\cdots\boxplus(\tau_r,2)
$$
where $\phi_\tau=(\tau_1,1)\boxplus\cdots\boxplus(\tau_r,1)$ is a generic global Arthur parameter of $G_n^*$. Then $\phi_\tau^{(2)}$ is a global Arthur parameter of $H_{2n}^*$. 
For any $\Sigma\in\CA_2(H_{2n})$ belonging to the global Arthur packet $\wt{\Pi}_{\phi_\tau^{(2)}}(H_{2n})$, 
a Bessel module $\CF^{\CO_\kappa}(\CC_\Sigma)$ with $\kappa=n$ or $n-1$, depending on the type of $H_{2n}^*$, 
is cuspidal on $H_{2n}^{\CO_\kappa}\cong G_n^*$. One of the technical difficulties in the theory is to prove that the Bessel module $\CF^{\CO_\kappa}(\CC_\Sigma)$ is in fact non-zero. This global non-vanishing property 
was established if one takes $\Sigma$ to be the residual representation $\CE_\tau$ in $\wt{\Pi}_{\phi_\tau^{(2)}}(H_{2n})$. In this case, the Bessel module $\CF^{\CO_\kappa}(\CC_{\CE_\tau})$ turns out to be 
an irreducible generic cuspidal representation of $G_n^*(\BA)$ with its global Arthur parameter $\phi_\tau$. 

The theory of twisted automorphic descents is to take automorphic representations in the global Arthur packet $\wt{\Pi}_{\psi^{(2,1)}}(H_{2n+m})$ as the source representations, where 
$$
\psi^{(2,1)}=\phi_\tau^{(2)}\boxplus\phi_\delta,
$$
and $\phi_\delta=(\delta_1,1)\boxplus\cdots\boxplus(\delta_s,1)$ is a generic global Arthur parameters of $H_m$. The key point in the theory of twisted automorphic descents is the input of the generic global Arthur 
parameter $\phi_\delta$. We fix $\phi_\tau^{(2)}$-part and let $\phi_\delta$ vary. The theory predicts the construction of all cuspidal members in the global Arthur packet $\wt{\Pi}_{\phi_\tau}(G_n)$ for 
pure inner forms $G_n$ of $G_n^*$. The proof of the global non-vanishing of the descent $\CF^{\CO_\ell}(\CC_{\CE_{\tau\otimes\sigma}})$ for $\sigma$ varying in the global Arthur packet 
$\wt{\Pi}_{\phi_\delta}(H_m)$ is directly related to the non-vanishing of the central value of the tensor product $L$-function $L(s,\tau\times\sigma)$. Since $\tau$ is fixed, but $\sigma$ varying, this $L$-function 
can also be called the standard $L$-function of $\tau$ twisted by $\sigma$. This is where the name of twisted automorphic descents comes from. We refer to \cite{JZ-BF}, \cite{J-Shahidi}, and \cite{JZ-Howe} for further discussion of the theory, and its applications.


\section{Branching Problem: Cuspidal Case}\label{sec-BPC}


For  $\pi\in\CA_\cusp(G_n)$ with cuspidal realization $\CC_\pi$, the {\sl Branching Problem} is to consider the spectral structure of the Bessel module $\CF^{\CO_\ell}(\CC_\pi)$ as a representation of 
$G_{\ell^-}^{\CO_\ell}(\BA)$. We first recall Proposition 2.2 of \cite{JZ-BF}, which establishes the tower property and the cuspidality of 
Bessel modules $\CF^{\CO_\ell}(\CC_\pi)$ of $\pi$. 

\begin{prop}[Tower Property and Cuspidality]\label{tpc}
Assume that $\pi\in\CA_\cusp(G_n)$ has a cuspidal realization $\CC_\pi$.  
Then the $\ell$-th Bessel module $\CF^{\CO_\ell}(\CC_\pi)$ of $\pi$ enjoys the following property:
There exists an integer $\ell_0$ in $\{0,1,\cdots,\Fr\}$, where $\Fr$ is the $k$-rank of $G_n$, such that
\begin{enumerate}
\item  the $\ell_0$-th Bessel module
$\CF^{\CO_{\ell_0}}(\CC_\pi)$ of $\CC_\pi$ is nonzero, but for any $\ell\in\{0,1,\cdots,\Fr\}$ with $\ell>\ell_0$, the $\ell$-th Bessel module
$\CF^{\CO_\ell}(\CC_\pi)$ is identically zero; and
\item the  $\ell_0$-th Bessel module $\CF^{\CO_{\ell_0}}(\CC_\pi)$ is cuspidal in the sense that its constant terms along all the parabolic subgroups of
$H_{\ell_0^-}^{\CO_{\ell_0}}$ are zero.
\end{enumerate}
\end{prop}

We call $\ell_0=\ell_0(\CC_\pi)$ in Proposition \ref{tpc} the {\sl first occurrence index} of $\CC_\pi$. When the realization $\CC_\pi$ is unique, we may write $\ell_0=\ell_0(\pi)$.
The tower property and the cuspidality of the $\ell$-th Bessel modules as in Proposition \ref{tpc} lead naturally to the following questions. 

\begin{ques}[First Occurrence Index]\label{q-foi}
For a given $\pi\in\CA_\cusp(G_n)$ with a cuspidal realization $\CC_\pi$, how to determine the first occurrence index $\ell_0=\ell_0(\CC_\pi)$ in terms of basic invariants of $\pi$?
\end{ques}

In terms of the global Arthur parameters of $\pi\in\CA_\cusp(G_n)$, Question \ref{q-foi} is closely related to Conjecture 4.2 of \cite{J14}. It will be very interesting to understand the 
relations between the first occurrence index $\ell_0=\ell_0(\CC_\pi)$ and the other basic invariants of $\pi$. Note that Conjecture 4.2 of \cite{J14} has been intensively studied by B. Liu and the first named author 
of this paper 
for symplectic groups and the metaplectic double cover of symplectic groups, but it is much less known for orthogonal groups. 

\begin{ques}[Spectrum]\label{q-sp}
For a given $\pi\in\CA_\cusp(G_n)$ with a cuspidal realization $\CC_\pi$, what can one say about the spectrum of the $\ell_0$-th Bessel module $\CF^{\CO_{\ell_0}}(\CC_\pi)$ as a representation of 
$H_{\ell_0^-}^{\CO_{\ell_0}}(\BA)$?
\end{ques}

This can be regarded as a cuspidal version of the branching problem (Question \ref{q-L2bp}) and may be quite difficult.
A preliminary understanding to Question \ref{q-sp} is the {\sl Generic Summand Conjecture} (Conjecture 2.3 of \cite{JZ-BF}). For completeness, we state it here. 

\begin{conj}[Generic Summand]\label{conj-gs}
Assume that $\pi\in\CA_\cusp(G_n)$ has a $G_n$-relevant, generic global Arthur parameter $\phi\in\wt{\Phi}_2(G_n^*)$. Then there exists a
cuspidal realization $\CC_\pi$ of $\pi$ in $L_\cusp^2(G_n)$ with the first occurrence index
$\ell_0=\ell_0(\CC_\pi)$, such that there exists a $k$-rational orbit
$\CO_{\ell_0}=\CO_{\udl{p}_{\ell_0}}$ in the $k$-stable orbits $\CO_{\udl{p}_{\ell_0}}^\st$ associated to the partition $\udl{p}_{\ell_0}$
with the {\bf Generic Summand Property:}
There exists at least one $\sigma$ in $\CA_\cusp(H_{\ell_0^-}^{\CO_{\ell_0}})$ with
an $H_{\ell_0^-}^{\CO_{\ell_0}}$-relevant, generic global Arthur parameter $\phi_\sigma$ in $\wt{\Phi}_2(H_{\ell_0^-}^*)$, and with a cuspidal
realization $\CC_\sig$ of $\sig$ in $L_\cusp^2(H_{\ell_0^-}^{\CO_{\ell_0}})$,
such that the $L^2$-inner product
$$
\left<\CF^{\psi_{\CO_{\ell_0}}}(\varphi_\pi),\varphi_\sigma\right>_{H_{\ell_0^-}^{\CO_{\ell_0}}}
$$
in the Hilbert space $L^2_\cusp(H_{\ell_0^-}^{\CO_{\ell_0}})$ is nonzero for some $\varphi_\pi\in\CC_\pi$ and $\varphi_\sigma\in\CC_\sigma$.
\end{conj}

It is well-known that one may explain the existence of such generic summand in the $\ell$-Bessel modules in terms of Vogan packets and of the non-vanishing of the central value of the 
tensor product $L$-function $L(s,\pi\times\sigma)$, according to the global Gan-Gross-Prasad conjecture. However, Conjecture \ref{conj-gs} asserts that if an irreducible cuspidal automorphic representation 
$\pi$ has a generic global Arthur parameter, then for one of its automorphic realization $\CC_\pi$, the Bessel module $\CF^{\CO_{\ell_0}}(\CC_\pi)$ at the first occurrence index $\ell_0=\ell_0(\CC_\pi)$, 
which is cuspidal, has at least one irreducible summand $\sigma$, which has a generic global Arthur parameter. The proof of Conjecture \ref{conj-gs} for general situation is still out of reach, 
although one may find the proof of some special cases in Section 7 of \cite{JZ-BF}. 

According to Theorem \ref{thm-psim}, with the assumption in Conjecture \ref{conj-gs}, the global Arthur parameter of any cuspidal summand $\sigma$ in the Bessel module $\CF^{\psi_{\CO_{\ell}}}(\CC_\pi)$ is of the form: 
$$
\psi_m=(\zeta_1,1)\boxplus\cdots\boxplus(\zeta_l,1)\boxplus(\xi_1,2)\boxplus\cdots\boxplus(\xi_k,2), 
$$
with $l,k\geq 0$. The reason, with which Conjecture \ref{conj-gs} should be true, is that  Conjecture \ref{conj-gs} only consider the spectrum of the 
Bessle module $\CF^{\psi_{\CO_{\ell}}}(\CC_\pi)$ at the first occurrence index $\ell_0$. This is compatible with Conjecture \ref{conj-pi-generic}.

One may consider a further refinement of the branching problem by requiring that both $\pi$ and $\sigma$ have generic global Arthur parameters. 
This leads to a close connection with the global Gan-Gross-Prasad conjecture. 

Assume that $\pi\in\CA_\cusp(G_n)$ has a generic global Arthur parameter $\phi_n\in\wt{\Phi}_2(G_n^*)$, which is $G_n$-relevant, and $\sigma\in\CA_\cusp(H_m)$, with $H_m=H_{\ell^-}^{\CO_{\ell}}$, has 
a generic global Arthur parameter $\phi_m\in\wt{\Phi}_2(H_m^*)$, which is $H_m$-relevant. The global Gan-Gross-Prasad conjecture implies that 
if the Bessel period $\CB^{\CO_{\ell}}(\varphi_\pi,\varphi_\sigma)$ in \eqref{bp-L2} is non-zero for some $\varphi_\pi\in\CC_\pi$ and $\varphi_\sigma\in\CC_\sigma$, 
then the central $L$-value $L(\frac{1}{2},\phi_G\times\phi_H)$ is non-zero. 
This is now a theorem (\cite{JZ-BF}). Moreover, in this case there exists one and at most one pair $(\pi,\sigma)$ in the global Vogan packet $\wt{\Pi}_{\phi_n\times\phi_m}[G_n^*\times H_m^*]$, 
such that the Bessel period $\CB^{\CO_{\ell}}(\varphi_\pi,\varphi_\sigma)$ is non-zero for some $\varphi_\pi\in\CC_\pi$ and $\varphi_\sigma\in\CC_\sigma$. Hence by means of the formulation of the 
global Gan-Gross-Prasad conjecture, the non-vanishing of the central $L$-value $L(\frac{1}{2},\phi_n\times\phi_m)$, plus the compatibility of the sign of the local symplectic root numbers 
for the pair $(\pi_v,\sigma_v)$ at all local places $v$ of $k$ 
will detect the automorphic branching property that the Bessel period $\CB^{\CO_{\ell}}(\varphi_\pi,\varphi_\sigma)$ is non-zero for some data with the unique choice of the pair $(\pi,\sigma)$ in the 
global Vogan packet $\wt{\Pi}_{\phi_n\times\phi_m}[G_n^*\times H_m^*]$. Hence in order to understand the branching decomposition problem in the framework of the global Gan-Gross-Prasad conjecture, 
it remains interesting to ask the following question. 

\begin{ques}[Generic Automorphic Branching Problem]\label{q-gabp}
For a $\pi\in\CA_\cusp(G_n)$ with a generic global Arthur parameter $\phi_n$ of $G_n^*$, which is $G_n$-relevant, and for an integer $\ell\in\{0,1,\cdots,\Fr\}$ with $m=\ell^-$, find a special orthogonal 
group $H_m$, such that $G_n\times H_m$ is relevant, and find a $\sigma\in\CA_\cusp(H_m)$ with a generic global Arthur parameter $\phi_m$ of $H_m^*$, which is $H_m$-relevant, such that 
\begin{enumerate}
\item the central $L$-value $L(\frac{1}{2},\phi_n\times\phi_m)$ is non-zero for certain choice of $\phi_m$; and 
\item the compatibility of the sign of the local symplectic root numbers for the pair $(\pi_v,\sigma_v)$ at all local places $v$ of the number field $k$. 
\end{enumerate}
\end{ques}

Following from the global Gan-Gross-Prasad conjecture, the representations $\sigma$ that satisfy the two conditions must occur in the automorphic branching decomposition of $\pi$, that is, the 
Bessel period $\CB^{\CO_{\ell}}(\varphi_\pi,\varphi_\sigma)$ is non-zero for some $\varphi_\pi\in\CC_\pi$ and $\varphi_\sigma\in\CC_\sigma$. It is a very interesting problem to construct explicitly such 
representations $\sigma$ via a variant of the twisted automorphic descent method. 



\section{Reciprocal Branching Problem}\label{sec-RBP}


Let $\sigma\in\CA_2(H_m)$ with an automorphic realization $\CC_\sig$. 
According to Question \ref{q-L2rbp}, the {\sl $L^2$-automorphic reciprocal branching problem} is to 
find a special orthogonal group $G_n$, such that $(G_n,H_m)$ is a relevant 
pair in the sense of global Gan-Gross-Prasad conjecture,  
and find $\pi\in\CA_2(G_n)$  such that the Bessel period $\CB^{\CO_{\ell}}(\varphi_\pi,\varphi_\sigma)$ is non-zero 
for some $\varphi_\pi\in\CC_\pi$ and $\varphi_\sigma\in\CC_\sigma$.
Note that when both $\sigma$ and $\pi$ are not cuspidal, the Bessel period $\CB^{\CO_{\ell}}(\varphi_\pi,\varphi_\sigma)$ needs to be defined with regularization. In this case, one may follow the 
compatibility of Bessel periods with parabolic induction as indicated in \cite{GJR04}, \cite{GJR05}, and \cite{GJR09}, for instance. Hence we may, without loss of generality, only consider the case 
when $\sigma$ is cuspidal and $\pi\in\CA_2(G_n)$.

\subsection{$\sigma$ has a generic global Arthur parameter}\label{ssec-sig-g}
We may consider first the case that $\sigma$ has a generic global Arthur parameter 
\begin{equation}\label{sig-gap}
\phi_\zeta=(\zeta_1,1)\boxplus\cdots\boxplus(\zeta_l,1)
\end{equation}
where $\zeta=\zeta_1\boxplus\cdots\boxplus\zeta_l$ is the generic isobaric sum of irreducible cuspidal automorphic representations $\zeta_j$ of $\GL_{\alpha_j}(\BA)$ with $\zeta_j\not\cong\zeta_i$ if $j\neq i$. 
By Theorem \ref{thm-psin}, if there is a special orthogonal group $G_n$ such that $(G_n,H_m)$ is relevant, and there is a $\pi\in\CA_\cusp(G_n)$ such that the Bessel period 
$\CB^{\CO_{\ell}}(\varphi_\pi,\varphi_\sigma)$ is non-zero for some $\varphi_\pi\in\CC_\pi$ and $\varphi_\sigma\in\CC_\sigma$, then the global Arthur parameter of $\pi$ must be of the form 
\begin{equation}\label{ap-psin1}
\psi_n=(\eta_1,1)\boxplus\cdots\boxplus(\eta_k,1)\boxplus(\delta_1,2)\boxplus\cdots\boxplus(\delta_l,2).
\end{equation}
We are going to show how to refine this structure via the construction in the theory of twisted automorphic descents. 

Take $\tau=\tau_1\boxplus\cdots\boxplus\tau_r$ be an isobaric automorphic representation of $\GL_a(\BA)$ with $\tau_i\not\cong\tau_j$ if $i\neq j$. Let $\GL_a\times H_m$ be a standard Levi subgroup of 
$H_{a+m}$. With $(\GL_a\times H_m,\tau\otimes\sig)$, we may define an Eisenstein series $E(\cdot,\phi_{\tau\otimes\sig},s)$. By \cite{JLZ13}, this Eisenstein series has a pole at $s=\frac{1}{2}$ if and only if 
the $L$-function $L(s,\tau,\rho)$ has a pole at $s=1$ and the central value $L(\frac{1}{2},\tau\otimes\sig)$ is non-zero. The condition that $L(s,\tau,\rho)$ has a pole at $s=1$ is equivalent to say that 
$\tau$ generates a generic global Arthur parameter for an special orthogonal group of the same type as $G_n$. This implies in particular that $a$ must be even. 
We may assume that $a\geq 2m$ if $H_m$ is an even special orthogonal group; and $a>2m$ if $H_m$ is an odd special orthogonal group, in the following discussion. We also assume that 
the central value $L(\frac{1}{2},\tau\otimes\sig)$ is non-zero. Under these two assumptions, we obtain that the Eisenstein series $E(\cdot,\phi_{\tau\otimes\sig},s)$ has a non-zero residue 
$\CE_{\tau\otimes\sig}$, which is square-integrable and has a global Arthur parameter of the form 
$
\phi_\tau^{(2)}\boxplus\phi_\zeta.
$

We may consider the global Arthur packet $\wt{\Pi}_{\phi_\tau^{(2)}\boxplus\phi_\zeta}(H_{a+m})$. The automorphic members $\Sigma$ in $\wt{\Pi}_{\phi_\tau^{(2)}\boxplus\phi_\zeta}(H_{a+m})$ could be cuspidal 
or residual representations. We may define $\CA_{\phi_\tau^{(2)}\boxplus\phi_\zeta}(H_{a+m})$ to be the set of all $\Sigma\in\CA_2(H_{a+m})$ that also belongs to $\wt{\Pi}_{\phi_\tau^{(2)}\boxplus\phi_\zeta}(H_{a+m})$.
The tower of the Bessel modules $\CB^{\CO_\kappa}(\Sigma)$ for all $\Sigma\in \CA_{\phi_\tau^{(2)}\boxplus\phi_\zeta}(H_{a+m})$ has been studied in \cite[Section 6.1]{JZ-BF}. 
In particular, Proposition 6.2 of \cite{JZ-BF} proves that for every $\Sigma\in \CA_{\phi_\tau^{(2)}\boxplus\phi_\zeta}(H_{a+m})$, 
the Bessel modules $\CB^{\CO_\kappa}(\Sigma)$ are zero if $\kappa>\frac{a}{2}+m-1$ when $H_m$ is an even special orthogonal group; and if $\kappa>\frac{a}{2}+m$ when $H_m$ is an odd special orthogonal group, respectively. 
However, it remains mysterious about 
the first occurrence index $\ell_0(\CC_\Sigma)$ of each member $\Sigma$ in $\CA_{\phi_\tau^{(2)}\boxplus\phi_\zeta}(H_{a+m})$. 
In the spirit of the general conjecture (\cite[Conjecture 4.2]{J14}), or more specifically, Conjecture 6.7 of \cite{JZ-BF}, there should be at least one member in $\CA_{\phi_\tau^{(2)}\boxplus\phi_\zeta}(H_{a+m})$, such that 
the first occurrence index $\ell_0(\CC_\Sigma)=\frac{a}{2}+m-1$ when $H_m$ is an even special orthogonal group; and $\ell_0(\CC_\Sigma)=\frac{a}{2}+m$ when $H_m$ is an odd special orthogonal group, respectively. 
This expectation has been verified in Proposition 6.8 of \cite{JZ-BF} for $m=0$ and $m=1$. The close relation of this expectation to the global Gan-Gross-Prasad conjecture has been explained in Theorem 6.9 and its proof 
in \cite{JZ-BF}. 

We define $\ell_0(\phi_\tau^{(2)}\boxplus\phi_\zeta)$ to be the largest possible first occurrence index of all 
$\Sigma$ in $\CA_{\phi_\tau^{(2)}\boxplus\phi_\zeta}(H_{a+m})$ and call it the first occurrence index of the global Arthur packet $\wt{\Pi}_{\phi_\tau^{(2)}\boxplus\phi_\zeta}(H_{a+m})$, 
which means that $\ell_0(\phi_\tau^{(2)}\boxplus\phi_\zeta)=\frac{a}{2}+m-1$ when $H_m$ is an even special orthogonal group; and 
$\ell_0(\phi_\tau^{(2)}\boxplus\phi_\zeta)=\frac{a}{2}+m$ when $H_m$ is an odd special orthogonal group. 
For $\Sigma=\CE_{\tau\otimes\sig}$ and $\sig\in \CA_{\phi_\zeta}(H_m)$, which belongs to 
$\CA_{\phi_\tau^{(2)}\boxplus\phi_\zeta}(H_{a+m})$, we consider the Bessel modules $\CF^{\CO_\kappa}(\CE_{\tau\otimes\sig})$ with $0\leq\kappa\leq \ell_0(\phi_\tau^{(2)}\boxplus\phi_\zeta)$. Assume that 
$\kappa_0=\ell_0(\CE_{\tau\otimes\sig})$ is the first occurrence index of the residual representation $\CE_{\tau\otimes\sig}$. Then by the same argument of the proof of Proposition 6.3 of \cite{JZ-BF}, the 
Bessel module $\CF^{\CO_{\kappa_0}}(\CE_{\tau\otimes\sig})$ is non-zero and cuspidal, and hence can be written as 
\begin{equation}\label{1stocc}
\CF^{\CO_{\kappa_0}}(\CE_{\tau\otimes\sig})
=
\pi_1\oplus\pi_2\oplus\cdots\oplus\pi_j\oplus\cdots
\end{equation}
as a representation of $G_{\kappa_0^-}^{\CO_{\kappa_0}}(\BA)$, where every $\pi_j$ is irreducible cuspidal automorphic representations of $G_{\kappa_0^-}^{\CO_{\kappa_0}}(\BA)$. Since the residual representation 
$\CE_{\tau\otimes\sig}$ is irreducible, by the uniqueness of local Bessel models at all local places of $k$, the decomposition in \eqref{1stocc} is of multiplicity free. It is clear that for each $\pi_i$, 
the Bessel period $\CB^{\CO_{\kappa_0}}(\CE_{\tau\otimes\sig},\varphi_{\pi_i})$ is non-zero for some choice of data. By the reciprocal non-vanishing of Bessel 
periods as established in Theorem 5.3 of \cite{JZ-BF}, there exists an integer $\ell_0$ such that the Bessel period $\CB^{\CO_{\ell_0}}(\varphi_{\pi_i},\varphi_\sig)$ is non-zero. This implies the following 
theorem, which can be viewed as an answer to the 
{\sl Cuspidal Reciprocal Branching Problem}. 

\begin{thm}[Cuspidal Reciprocal Branching]\label{thm-crb}
Assume that $\sig\in\CA_\cusp(H_m)$ has a generic global Arthur parameter $\phi_\delta$ as in \eqref{sig-gap}. 
For any even integer $a$, which is greater than or equal to $2m$ if $H_m$ is an even special orthogonal group; and is greater than $2m$ if $H_m$ is an odd special orthogonal group, assume that there exists 
an isobaric automorphic representation $\tau=\tau_1\boxplus\cdots\boxplus\tau_r$ of $\GL_a(\BA)$ with the properties that the global Arthur parameter $\phi_\tau$ is of the type opposite to that of $\phi_\delta$, 
and the central value $L(\frac{1}{2},\tau\times\sig)$ is not zero. 
Then there exists a special 
orthogonal group $G_n$ such that $(G_n,H_m)$ is a relevant pair, and there is a $\pi\in\CA_\cusp(G_n)$ with an automorphic realization $\CC_\pi$, whose global Arthur parameter $\psi_n$ 
is of the form as in \eqref{ap-psin1}, such that 
the Bessel period $\CB^{\CO_{\ell_0}}(\varphi_{\pi},\varphi_\sig)$ is non-zero for some choice of $\varphi_\sig\in\CC_\sig$ and $\varphi_\pi\in\CC_\pi$.
\end{thm}

Under the assumption, the residue $\CE_{\tau\otimes\sig}$ is non-zero and belongs to the global Arthur packet $\wt{\Pi}_{\phi_\tau^{(2)}\boxplus\phi_\delta}(H_{a+m})$. Assume that 
$\kappa_0=\ell_0(\CE_{\tau\otimes\sig})$ is the first occurrence index of the residual representation $\CE_{\tau\otimes\sig}$. Then $\Fn=2a+\Fm-2\kappa_0-1$, $n=[\frac{\Fn}{2}]$, and 
$\ell_0=a-\kappa_0-1$. 

Let us discuss a special case when the first occurrence index $\kappa_0=\ell_0(\CE_{\tau\otimes\sig})$ is equal to $\ell_0(\phi_\tau^{(2)}\boxplus\phi_\delta)$, 
the first occurrence index of the global Arthur packet $\wt{\Pi}_{\phi_\tau^{(2)}\boxplus\phi_\zeta}(H_{a+m})$, then $a=2n$ and $\kappa_0^-=n$. In this case, by Theorem 6.5 of \cite{JZ-BF}, the irreducible summands $\pi_i$ of 
$\CF^{\CO_{\kappa_0}}(\CE_{\tau\otimes\sig})$ in \eqref{1stocc} belong to $\CA_\cusp(G_{n}^{\CO_{\kappa_0}})$ have $\phi_\tau$ as their global Arthur parameter. Again, by the reciprocal non-vanishing 
of Bessel Periods (Theorem 5.3 of \cite{JZ-BF}), the Bessel periods $\CB^{\CO_{\ell_0}}(\varphi_{\pi_i},\varphi_\sig)$ is non-zero for some choice of data. Assuming the local Gan-Gross-Prasad conjecture 
for special orthogonal groups over all local fields, we deduce that $\CF^{\CO_{\kappa_0}}(\CE_{\tau\otimes\sig})=\pi$ is irreducible if $G_{n}^{\CO_{\kappa_0}}$ is an odd special orthogonal group; and 
$\CF^{\CO_{\kappa_0}}(\CE_{\tau\otimes\sig})=\pi\oplus\pi^\star$ if $G_{n}^{\CO_{\kappa_0}}$ is an even special orthogonal group. The details of the proof for this last assertion is the same as the proof of 
Theorem 7.1 of \cite{JZ-BF}. 

\begin{thm}\label{thm-crb-foc}
Take the same assumption as in Theorem \ref{thm-crb}. If the first occurrence index $\kappa_0=\ell_0(\CE_{\tau\otimes\sig})$ is equal to $\ell_0(\phi_\tau^{(2)}\boxplus\phi_\delta)$, 
the first occurrence index of the global Arthur packet $\wt{\Pi}_{\phi_\tau^{(2)}\boxplus\phi_\zeta}(H_{a+m})$, then there exists a $\pi\in\CA_\cusp(G_{n}^{\CO_{\kappa_0}})$ with $2n=a$ with $\phi_\tau$ as its generic global Arthur parameter, such that the Bessel period $\CB^{\CO_{\ell_0}}(\varphi_{\pi},\varphi_\sig)$ is non-zero for some choice of $\varphi_\sig\in\CC_\sig$ and $\varphi_\pi\in\CC_\pi$.
\end{thm}

\subsection{Relations with the global Gan-Gross-Prasad conjecture}\label{ssec-rggp}
We are going to discuss the relation between the {\sl Cuspidal Reciprocal Branching Problem} and the global Gan-Gross-Prasad conjecture. 

Given $\sigma\in\CA_\cusp(H_m)$ with an $H_m$-relevant, generic global Arthur parameter $\phi_\delta$ of $H_m^*$, as in \eqref{sig-gap}, 
it is not hard to find a family of special orthogonal groups $G_n$ such that $(G_n,H_m)$ is a relevant pair in the sense of global Gan-Gross-Prasad conjecture. It is itself a hard problem to find 
a $G_n$-relevant, generic global Arthur parameter $\phi_\tau$ of $G_n^*$ with $a=2n$, such that 
the central $L$-value $L(\frac{1}{2},\phi_\tau\times\phi_\delta)$ is non-zero. We take it as an assumption. 
The global Gan-Gross-Prasad conjecture asserts that there exists a unique pair 
$(\pi_0,\sigma_0)$ in the global Vogan packet $\wt{\Pi}_{\phi_G\times\phi_H}[G_n^*\times H_m^*]$ with the property that the Bessel period 
$\CB^{\CO_{\ell_0}}(\varphi_{\pi_0},\varphi_{\sigma_0})$, with $m=\ell_0^-$, 
is non-zero for some $\varphi_{\pi_0}\in\CC_{\pi_0}$ and $\varphi_{\sigma_0}\in\CC_{\sigma_0}$.

If the given $\sigma$, sharing the global Vogan packet with $\sigma_0$, is happened to be $\sigma_0$, 
then the cuspidal automorphic module $\pi$ as constructed by the twisted automorphic descent in 
Theorem \ref{thm-crb-foc} is equal to $\pi_0$. 

If $\sigma$ shares a global Vogan packet with $\sigma_0$, but they are not equivalent to each other, the global Gan-Gross-Prasad conjecture asserts that the first occurrence index 
$\kappa_0=\ell_0(\CE_{\tau\otimes\sig})$ must be  less than $\ell_0(\phi_\tau^{(2)}\boxplus\phi_\delta)$. This goes to the situation considered in Theorem \ref{thm-crb}. In this case, we may still apply the 
twisted automorphic descent method to construct an irreducible cuspidal automorphic representation $\pi$ of $G_n$, such that 
the Bessel period $\CB^{\CO_{\ell_0}}(\varphi_{\pi},\varphi_\sig)$ is non-zero for some choice of $\varphi_\sig\in\CC_\sig$ and $\varphi_\pi\in\CC_\pi$. However, we do not have precise information about the 
global Arthur parameter $\psi_n$ of $\pi$, except that $\psi_n$ is of the form in \eqref{ap-psin1}. Hence it is a desirable question to find the precise information about the global Arthur parameter of $\pi$ that gives 
an answer the {\sl Cuspidal Reciprocal Branching Problem}!
In a recent work of Jiang-Liu-Xu (\cite{JLX}), a special situation when $H_m$ is $\SO(3)$ is explicitly studied. 
In this case, we are able to prove in some situation that the twisted automorphic descent method may still 
construct a cuspidal $\pi$ with generic global Arthur parameter (Theorem 6.3 of \cite{JLX}). We will leave further details of discussion about this mysterious question to our future work. 

On the other hand, Theorem \ref{thm-crb} can be viewed as examples of the global Gan-Gross-Prasad conjecture for non-generic global Arthur parameters, although the correct formulation of the new conjecture is 
still not known (\cite{GGP18}).


\end{document}